%% file: 2005-50.tex
\documentclass{gtart_h}  

\input gtoutput

\lognumber{528}
\received{25 November 2004}
\volumenumber{9}\papernumber{50}\volumeyear{2005}
\pagenumbers{2193}{2226}   
\revised{4 September 2005}
\published{1 December 2005}
\accepted{26 November 2005}
\proposed{Yasha Eliashberg}
\seconded{Tomasz Mrowka, Joan Birman}

\usepackage{epsf,amssymb,amsmath}    

\newtheorem{theoreme}{Th\'eor\`eme}[section]
\newtheorem{proposition}[theoreme]{Proposition}

\newtheorem{lemme}[theoreme]{Lemme}
\newtheorem{corollaire}[theoreme]{Corollaire}
\newtheorem{conjecture}[theoreme]{Conjecture}
\newtheorem{question}[theoreme]{Question}
\newtheorem{remarque}[theoreme]{Remarque}
\newtheorem{affirmation}[theoreme]{Affirmation}

\numberwithin{equation}{subsection}

\newcommand{\R}{\mathbb{R}}
\newcommand{\Z}{\mathbb{Z}}
\newcommand{\Q}{\mathbb{Q}}
\newcommand{\N}{\mathbb{N}}

\newcommand{\bdry}{\partial}
\newcommand{\s}{\medskip}
\newcommand{\n}{\noindent}
\newcommand{\B}{\mathcal{B}}
\newcommand{\F}{\mathcal{F}}

\newcommand{\be}{\begin{enumerate}}
\newcommand{\ee}{\end{enumerate}}

\begin{document}

\title{Constructions contr\^ol\'ees de champs de Reeb\\et applications}
\covertitle{Constructions contr\noexpand\^ol\noexpand\'ees de 
champs de Reeb\\et applications}
\asciititle{Constructions controlees de champs de Reeb et applications}

\author{Vincent Colin\\Ko Honda}

\address{Universit\'e de Nantes, UMR 6629 du CNRS, 44322 Nantes,
  France\\
University of Southern California, Los Angeles, CA 90089, USA}
\asciiaddress{Universite de Nantes, UMR 6629 du CNRS, 44322 Nantes,
  France\\and\\University of Southern California, Los Angeles, CA
  90089, USA}

\gtemail{\mailto{Vincent.Colin@math.univ-nantes.fr}{\rm\qua 
and\qua}\mailto{khonda@math.usc.edu}}
\asciiemail{Vincent.Colin@math.univ-nantes.fr, khonda@math.usc.edu}        
\urladdr{http://rcf.usc.edu/~khonda}



\begin{abstract}\vspace{-1mm}
On every compact, orientable, irreducible 3--manifold $V$ which is
toroidal or has torus boundary components we construct a contact 1--form
whose Reeb vector field $R$ does not have any contractible periodic
orbits and is tangent to the boundary.  Moreover, if $\bdry V$ is
nonempty, then the Reeb vector field $R$ is transverse to a taut
foliation.  By appealing to results of Hofer, Wysocki, and Zehnder, we
show that, under certain conditions, the 3--manifold obtained by Dehn
filling along $\bdry V$ is irreducible and different from the 3--sphere.

{\bf R\'esum\'e}

On construit, sur toute vari\'et\'e $V$ de dimension trois orientable, 
compacte, irr\'eductible, bord\'ee par des tores 
ou toro\"\i dale, une
forme de contact dont le champ de Reeb $R$ est sans orbite p\'eriodique
contractible et tangent au bord.
De plus, si $\partial V$ est non vide, le champ $R$
est transversal \`a un feuilletage tendu.
En utilisant des r\'esultats de Hofer, Wysocki et Zehnder, on obtient
sous certaines conditions que la vari\'et\'e obtenue par
obturation de Dehn le long du bord de $V$ est irr\'eductible et
diff\'erente de la sph\`ere $S^3$.
\end{abstract}     

\asciiabstract{%
On every compact, orientable, irreducible 3-manifold V which is
toroidal or has torus boundary components we construct a contact 1-form
whose Reeb vector field R does not have any contractible periodic
orbits and is tangent to the boundary.  Moreover, if bdry V is
nonempty, then the Reeb vector field R is transverse to a taut
foliation.  By appealing to results of Hofer, Wysocki, and Zehnder, we
show that, under certain conditions, the 3-manifold obtained by Dehn
filling along bdry V is irreducible and different from the 3-sphere.

Resume

On construit, sur toute variete V de dimension trois orientable,
compacte, irreductible, bordee par des tores ou toroidale, une forme
de contact dont le champ de Reeb R est sans orbite periodique
contractible et tangent au bord.  De plus, si le bord de V est non
vide, le champ R est transversal a un feuilletage tendu.  En utilisant
des resultats de Hofer, Wysocki et Zehnder, on obtient sous certaines
conditions que la variete obtenue par obturation de Dehn le long du
bord de V est irreductible et differente de la sphere S^3.}

\maketitle

\section{Introduction}

La topologie de contact a connu ces derni\`eres ann\'ees des 
\'evolutions spectaculaires. L'objet de ce texte est de conforter les
liens qui la relient \`a la topologie de dimension trois.

Plus pr\'ecisemment, on inaugure l'\'etude des structures
de contact {\it hypertendues}, dans une tentative de 
rapprocher deux facettes de la g\'eom\'etrie de contact: la dynamique
et la th\'eorie des feuilletages.

On se place sur une vari\'et\'e de dimension trois.

Une structure de contact dont un champ de Reeb est
sans orbite p\'eriodique contractible sera dite {\it hypertendue}.
Lorsqu'on veut pr\'eciser le champ de Reeb,
on parle de forme de contact hypertendue.
Toute structure de contact hypertendue est tendue d'apr\`es le
th\'eor\`eme fondamental de H~Hofer, K~Wysocki et E~Zehnder:

\begin{theoreme}\label{theoreme : Hofer} {\rm \cite{Hof,HWZ1,HWZ2}}    

{\rm(a)}\qua Si $R$ est un champ de Reeb sans orbite p\'eriodique contractible
sur une vari\'et\'e de contact $(V,\xi )$ close de dimension trois, alors $V$ est
irr\'eductible et $\xi$ est tendue.

{\rm(b)}\qua Tout champ de Reeb $R$ sur la sph\`ere
$S^3$ poss\`ede une orbite p\'eriodique non nou\'ee.
De plus, si le champ $R$ est non d\'eg\'en\'er\'e,
il existe une orbite p\'eriodique qui borde un disque plong\'e dont
l'int\'erieur est transversal \`a $R$.
Cette propri\'et\'e est \'egalement v\'erifi\'ee
pour les formes de contact vrill\'ees g\'en\'eriques
sur les sph\`eres d'homologie.

\end{theoreme}

On peut lire la partie (b) de ce r\'esultat \`a l'envers:
si $R$ est un champ de Reeb g\'en\'erique sur $V$ dont aucune 
orbite p\'eriodique contractible ne borde de disque immerg\'e $D$,
d'int\'erieur transversal \`a $R$, alors
le rev\^etement universel de $V$ n'est pas $S^3$.

Ajout\'e au crit\`ere d'irr\'eductibilit\'e,
on peut faire la conjecture suivante:

\begin{conjecture} \label{univcover}
Si $V$ porte une forme de contact hypertendue,
alors le rev\^ete\-ment universel de $V$ est $\R^3$.
\end{conjecture}                   

Au regard des travaux r\'ecents de Perelman \cite{Pe1,Pe2,Pe3}, 
la
Conjecture~\ref{univcover} pourrait en fait \^etre d\'ej\`a \'etablie.

Dans ce texte, les structures de contact hypertendues sont appel\'ees
\`a jouer un r\^ole similaire \`a celui des
feuilletages (de codimension $1$) tendus dans les travaux de 
D~Gabai~\cite{Ga1,Ga2,Ga3}.
La d\'efinition de structure
hypertendue fait en particulier \'echo \`a un r\'esultat fondamental de
S~Novikov \cite{Nov}, qui implique que toute courbe transversale
\`a un feuilletage tendu est non contractible.
Y Eliashberg et W Thurston \cite{ET} ont par ailleurs
montr\'e que toute structure de contact $C^0$--proche
d'un feuilletage tendu est universellement tendue.   

Notre travail est aussi motiv\'e par des travaux r\'ecents de 
D Gabai--L Mosher, S Fenley, W Thurston et D Calegari (voir
\cite{Mo,Fe,Ca1,Ca2,Ca3, Th2}), dans lesquels ces auteurs \'etudient
les flots pseudo-Anosov qui sont transversaux (ou presque transversaux)
\`a un feuilletage tendu.
Dans notre situation, l'\'etude topologique des structures de
contact tendues est compl\'et\'ee par l'\'etude dynamique des
champs de Reeb.

Ici, on montre l'existence de structures de contact hypertendues.
En comparant notre construction avec celle de Gabai~\cite{Ga1}, 
on obtient \'egalement un contr\^ole suppl\'ementaire sur le champ de 
Reeb: il est transversal \`a un feuilletage tendu.

\begin{theoreme}\label{theoreme : existence} 
Toute vari\'et\'e $V$ compacte, orientable,
irr\'eductible de dimension $3$, bord\'ee par une union non vide
de tores, porte une forme de contact hypertendue,
dont le champ de Reeb est tangent au bord et transversal \`a un feuilletage
tendu $\F$.
\end{theoreme}

\begin{question}
{\rm Sur une vari\'et\'e close et orientable, tout feuilletage de 
codimension $1$ sans feuille compacte est-il transversal \`a un 
champ de Reeb?}
\end{question}     

Plus vraisemblablement, on peut envisager que:

\begin{conjecture}
{\rm Toute vari\'et\'e close et orientable qui porte un feuilletage 
tendu, porte aussi une structure de contact hypertendue.}
\end{conjecture}

La preuve du th\'eor\`eme~\ref{theoreme : existence} repose sur 
l'existence d'une hi\'erarchie sutur\'ee de $V$, comme dans \cite{Ga1}.
Comme corollaire facile, on obtient:

\begin{corollaire}\label{corollaire : toroidal} 
Toute vari\'et\'e close, orientable, irr\'eductible et toro\" \i dale 
de dimension $3$ porte une structure de contact hypertendue.
\end{corollaire}                  

On retrouve ainsi un r\'esultat de \cite{Co3} (voir \cite{HKM2}
pour une autre d\'emonstration).

Les champs de Reeb produits par le th\'eor\`eme~\ref{theoreme : existence}
ont une dynamique particuli\`erement bien contr\^ol\'ee qui les
rend adapt\'es au calcul de leur {\em homologie de contact}~\cite{EGH}.
Ce projet est men\'e \`a bien dans un article en pr\'eparation
du premier auteur avec F~Bourgeois~\cite{BC}, qui montre que l'homologie de
contact distingue une infinit\'e de structures de contact hypertendues
sur toute vari\'et\'e orientable, close, irr\'eductible et toro\"\i dale de 
dimension trois. (Ce calcul donne une autre preuve du fait, d\^u \`a
Colin~\cite{Co5} et Honda--Kazez--Mati\'c~\cite{HKM2}, que les vari\'et\'es
toro\"\i dales portent une infinit\'e de structures de contact tendues.)  
On d\'emontre ainsi
 \'egalement la conjecture de Weinstein pour presque toutes les structures de contact
de torsion non nulle connues sur les vari\'et\'es toro\"\i dales.

On applique ces id\'ees sur un exemple. Soit $S$ une surface
compacte orientable \`a bord non vide, de caract\'eristique $\chi
(S)<0$, et $V$ la fibration sur le cercle obtenue par suspension
d'un diff\'eomorphisme de $S$, qui est l'identit\'e pr\`es de
$\partial S$. On note $T_1$,...,$T_n$ les composantes de $\partial V$.  
Pour tout $1\leq j\leq n$, on note $[m_j ]$ et $[l_j ]$ une base 
de $H_1 (T_j ;\Z )$, avec $i([m_j ],[l_j] )=1$ (pour l'orientation de 
$T_j$ comme bord de $V$), et $[m_j] =[\partial S \cap T_j ]$. Si 
$\gamma$ est une courbe non contractible dans $T_j$, sa pente est 
d\'efinie comme $p_j /q_j$, o\`u $[\gamma ]=p_j [l_j ] +q_j [m_j ]$. En
particulier, la pente de $\partial S \cap T_j$ vaut $0$.
Pour $\varepsilon =(\varepsilon_1 ,...,\varepsilon_n )\in \Q^n$, on note
$V_\varepsilon$ la vari\'et\'e  obtenue par obturation de Dehn le
long de $T_j$ avec la pente $\varepsilon_j$.

\begin{theoreme}\label{theoreme : obturation}
Si $\rho$ est assez proche de $0$ (pas forc\'ement positif),
et si pour tout $1\leq i\leq n$, $\varepsilon_i$ est non nul et compris
entre $0$ et $\rho$, alors $V_\varepsilon$ porte une structure de
contact hypertendue (positive ou n\'egative suivant le signe de 
$\rho$). La vari\'et\'e $V_\varepsilon$ est donc irr\'eductible et 
diff\'erente de $S^3$.
\end{theoreme}                                        

En particulier, pour $\varepsilon$ proche de $0$ et non nul, une
$\varepsilon$--chirurgie de Dehn sur un n\oe ud (et m\^eme un
entrelacs) fibr\'e non trivial de $S^3$ n'est pas $S^3$. Cette
application peut-\^etre vue comme le pendant, par le biais des
structures de contact, de la preuve par D~Gabai~\cite{Ga2} que la
chirugie d'indice $0$ sur un n\oe ud non trivial de $S^3$ ne donne
jamais $S^1 \times S^2$. Le th\'eor\`eme~\ref{theoreme :
obturation} donne de nombreux exemples de structures hypertendues
sur des sph\`eres d'homologie atoro\"\i dales et des vari\'et\'es
hyperboliques.   

Pour d\'emontrer que la forme de contact $\alpha_\varepsilon$
produite dans le th\'eor\`eme~\ref{theoreme : obturation} est
hypertendue, on utilise l'existence d'un feuilletage tendu sur
$V_\varepsilon$, construit par R.\ Roberts dans \cite{Ro1},
transversal au champ de Reeb $R_\varepsilon =R_{\alpha_\varepsilon}$.
Les cons\'equences topologiques dues \`a la pr\'esence de ce
feuilletage tendu recouvrent celles dues \`a la pr\'esence d'une
structure hypertendue: en fait, $V_\varepsilon$ est rev\^etue par
$\R^3$.

Une motivation pour utiliser la g\'eom\'etrie de contact vient du
fait que, d'apr\`es Hofer, Wysocki et Zehnder, la sph\`ere $S^3$
est compl\`etement caract\'eris\'ee par les propri\'et\'es
dynamiques des champs de Reeb qu'elle supporte. Par ailleurs,
l'extension de la structure de contact construite sur $V$
dans le th\'eor\`eme \ref{theoreme : existence}  en une structure
``potentiellement" (hyper-)tendue sur $V_\varepsilon$ est
extr\^emement simple, si bien qu'on peut esp\'erer appliquer ces
m\'ethodes de contact dans des situations o\`u les feuilletages
tendus sont inop\'erants.

En fait ici, on peut montrer directement, sans recours aux
feuilletages, que le champ de Reeb $R_\varepsilon$ produit sur
$V_\varepsilon$ ne poss\`ede pas d'orbite p\'eriodique qui borde un
disque plong\'e dont l'int\'erieur est transversal \`a
$R_\varepsilon$. Ceci suffit pour montrer, en utilisant la partie (b)
du th\'eor\`eme~\ref{theoreme : Hofer}, que pour tout $N \in \N$, si
$\rho$ est assez proche de $0$, la vari\'et\'e $V_\varepsilon$, ainsi
que tous ses rev\^etements de degr\'es inf\'erieurs \`a $N$, sont
diff\'erents de la sph\`ere $S^3$. La preuve, si elle est plus
satisfaisante pour qui cherche a \'eliminer les feuilletages et
obtenir des cons\'equences topologiques de l'\'etude des
structures de contact, est cependant assez technique et c'est
pourquoi on la reproduit seulement en annexe de ce texte. Il faut
donc voir l'utilisation des feuilletages comme un interm\'ediaire
pratique qui permet de r\'ev\`eler une propri\'et\'e du champ de
Reeb. En ce sens, on peut r\'eellement lire
l'\'enonc\'e~\ref{theoreme : obturation} comme une application
topologique de la g\'eom\'etrie de contact. On peut noter que le
disque d'int\'erieur transversal \`a $R$ et bord\'e par une orbite
p\'eriodique fourni par le th\'eor\`eme~\ref{theoreme : Hofer}
qu'on utilise dans l'annexe est en fait une feuille d'un feuilletage d'\'energie
finie de $S^3$ (voir \cite{HWZ2}), dont on pourrait tenter d'exploiter l'existence.

Pour en finir avec ces remarques, on note \'egalement qu'on
peut voir le th\'eor\`eme \ref{theoreme : obturation} comme
un cas tr\`es particulier d'une somme de r\'esultats dus \`a Thurston 
\cite{Th3} et \`a Culler, Gordon, Luecke et Shalen \cite{CGLS}:
si $V$ est atoro\"\i dale, elle est hyperbolique, et seules un 
nombre fini d'obturations ne sont pas hyperboliques et donc 
non rev\^etues par $\R^3$. Si $V$ est toro\"\i dale, seules un 
nombre fini d'obturations ne sont pas toro\"\i dales et donc
non rev\^etues par $\R^3$.

\begin{question}
{\rm Peut-on trouver une estimation pour la taille maximale de 
$\rho$, comme dans \cite{Ro2}?}
\end{question}

La structure $(V,\xi ,R,\F )$ fournie par le th\'eor\`eme
\ref{theoreme : existence} -- o\`u $R$ est un champ
de Reeb transversal \`a un feuilletage tendu $\F$ -- \'etendue
de mani\`ere standard \`a une obturation de Dehn de $V$ comme dans
le paragraphe~9, g\'en\'eralise
la notion de structure de contact port\'ee par un livre
ouvert d\'efinie par E~Giroux~\cite{Gi2}, o\`u le
feuilletage $\F$ est le feuilletage par les
fibres d'une fibration sur le cercle.
Elle rappelle \'egalement la notion de feuilletage d'\'energie
finie de Hofer, Wysocki et Zehnder~\cite{HWZ2}.

{\bf Remerciements}\qua VC remercie le CNRS pour son accueil
en d\'el\'egation au cours de l'ann\'ee 2003--2004. 
KH a \'et\'e subventionn\'e par la bourse DMS-023738
de la NSF et par la bourse Sloan.

\section{Sutures et hi\'erarchies sutur\'ees}
Une {\it vari\'et\'e sutur\'ee} $(V,\gamma )$ est la
donn\'ee d'une vari\'et\'e \`a coins (les coins sont model\'es sur
$(\R_+ )^2 \times \R$) 
de dimension $3$ compacte
et orient\'ee $V$ et d'un ensemble $\gamma \subset \partial V$
constitu\'e d'une collection d'anneaux deux \`a deux disjoints $A(\gamma )$
et de tores $T(\gamma )$.
Les coins de $V$ co\"\i ncident avec $\partial A(\gamma )$.

Chaque composante de $R(\gamma )=\partial V\setminus int (\gamma )$
est orient\'ee. On note $R_+ (\gamma )$ les composantes de $R(\gamma )$
le long desquelles le vecteur normal direct sort de $V$ et $R_- (\gamma )$ les
autres.

Chaque composante de $A(\gamma )$ est appel\'ee  {\it suture annulaire} et contient une {\it suture}, c'est-\`a-dire
une courbe orient\'ee, ferm\'ee simple et homologiquement non
triviale dans $A(\gamma )$. On note $s(\gamma )$ l'ensemble des sutures.

Les orientations de $R_+ (\gamma )$, $R_- (\gamma )$ et
$s(\gamma )$ v\'erifient la condition de compatibilit\'e
suivante: si $\alpha \subset \partial V$ est un arc orient\'e avec
$\partial \alpha \subset R(\gamma )$, qui a une intersection $+1$
avec $s(\gamma )$, alors $\alpha$ d\'ebute dans $R_- (\gamma )$
et aboutit dans $R_+ (\gamma )$.

On d\'efinit \`a pr\'esent le d\'ecoupage d'une vari\'et\'e
sutur\'ee $(V,\gamma )$ en $(V',\gamma ')$ le long d'une
surface $S$. La notation $N(A)$ d\'esigne un petit voisinage
tubulaire de $A$.

Soit $S$ une surface orient\'ee proprement plong\'ee dans $V$
avec les propri\'et\'es suivantes:
\begin{itemize}
\item $\partial S \pitchfork \gamma$;
\item si $S$ intersecte une suture annulaire $A\in A(\gamma )$
le long d'arcs, alors aucun d'eux ne s\'epare $A$;
\item si $S$ intersecte une suture annulaire $A$
le long de cercles, alors chacun d'eux, muni de l'orientation
d\'eduite de celle de $S$, est homologue \`a l'\^ame
$s(\gamma )\cap A$;
\item aucune composante de $S$ n'est un disque de bord
inclus dans $R(\gamma )$;
\item aucune composante de $\partial S$ ne borde de disque
dans $R(\gamma )$.
\end{itemize}

On note $V'$ la vari\'et\'e \`a bord
$V\setminus (int (N(S))$ et $S'_+$ et $S'_-$
les composantes de $\partial N(S)$ dans $\partial V'$,
o\`u le vecteur normal direct \`a $S$ pointe
respectivement, \`a l'ext\'erieur et \`a l'int\'erieur de $V'$.

On pose de plus
\begin{gather*} \gamma' =((\gamma \cap V') \cup N(S'_+ \cap R_- (\gamma ))
\cup N(S'_- \cap R_+ (\gamma )))\cap \partial V',\\
 R_+ (\gamma ')=((R_+ (\gamma )\cap V')\cup S'_+ ) -int (\gamma') ,\\
 R_- (\gamma ')=((R_- (\gamma )\cap V')\cup S'_- ) -int (\gamma') .
\end{gather*}
Une hi\'erarchie de vari\'et\'es sutur\'ees est une suite de tels
d\'ecoupages
$$(V,\gamma)=(V_0 ,\gamma_0 )\stackrel {S_0}\rightsquigarrow (V_1,\gamma_1)
\stackrel{S_1}\rightsquigarrow...\stackrel{S_{n-1}}\rightsquigarrow
(V_n ,\gamma_n )$$ le long de surfaces $S_0$,...,$S_n$, chaque
$S_i$ \'etant $\pi_1$--inject\'ee dans $V_i$, et aboutissant \`a
une union de boules $V_n$ dont les sutures $\gamma_n$ sur chaque
bord sont non vides et connexes (i.e.\ une collection de  $(D^2 \times [0,1],
\partial D^2 \times [0,1])$).

Une  vari\'et\'e sutur\'ee $(V,\gamma )$ est {\it annulaire},
si toutes les composantes de $V$ ont un bord non vide,
si toute composante de $\partial V$ contient une suture
et si toutes les sutures sont annulaires.
Elle est {\it tendue} si $V$ est irr\'eductible, $R(\gamma )$
est incompressible et minimise la norme de Thurston
dans $H_2 (V,\gamma;\Z )$.
Une hi\'erarchie sutur\'ee annulaire
$$(V,\gamma)=(V_0 ,\gamma_0 )\stackrel {S_0}\rightsquigarrow (V_1,\gamma_1)
\stackrel {S_1}\rightsquigarrow ...\stackrel {S_{n-1}}\rightsquigarrow 
(V_n ,\gamma_n )$$
est {\it bien positionn\'ee}
si toute composante de $\partial S_i$
rencontre $\gamma_i$ le long d'une famille non vide d'arcs.

\begin{theoreme}\label{theoreme : hierarchie} {\rm \cite{Ga1, HKM1}}\qua 
Soit $(V,\gamma )$ une vari\'et\'e sutur\'ee
annulaire tendue avec $H_2 (V,\partial V;\Z) \neq 0$. Alors
$(V,\gamma )$ admet une hi\'erarchie sutur\'ee annulaire tendue
bien positionn\'ee.
\end{theoreme}

Pour finir, on dit qu'une vari\'et\'e sutur\'ee $(V,\gamma )$ {\it porte}
un feuilletage $\F$ si $\F$ est d\'efini sur tout $V$,
transversal \`a $\gamma$, et si $R(\gamma )$ est une union de 
feuilles dont les orientations co\"\i ncident avec celles
de $\F$.

\section{Structures convexes et hi\'erarchies convexes}

La notion de hi\'erarchie sutur\'ee a \'et\'e adapt\'ee
aux vari\'et\'es de contact par K Honda, W Kazez et G Mati\'c
\cite{HKM1,HKM2}.
Une structure convexe sur une vari\'et\'e $V$ compacte orient\'ee
de dimension $3$ est la donn\'ee dans chaque composante
de $\partial V$ d'une collection non vide de courbes orient\'ees.
On note $\Gamma$ la r\'eunion de ces courbes, appel\'ees
courbes de s\'eparation.
Le compl\'ementaire $R=\partial V -\Gamma$ de $\Gamma$ dans
$\partial V$ est orient\'e.
On note $R_+$ la partie de $R$ o\`u le vecteur normal
\`a $R$ sort de $V$ et $R_- =R-R_+$.
Ces orientations v\'erifient la m\^eme condition
de compatibilit\'e que l'\^ame des sutures:
$\Gamma$ est le bord orient\'e de
l'adh\'erence de $R_+$ dans $\partial V$.

Soit $S$ une surface orient\'ee, proprement
plong\'ee dans $(V,\Gamma )$ et qui v\'erifie les propri\'et\'es suivantes:
\begin{itemize}
\item $\partial S \pitchfork \Gamma$;
\item aucune composante de $S$ n'est un disque de bord
inclus dans $R(\Gamma )$;
\item aucune composante de $\partial S$ ne borde de disque
dans $R(\Gamma )$.
\end{itemize}
On d\'efinit le d\'ecoupage $(V',\Gamma' )$ de $(V,\Gamma )$
le long de $S$ comme:
\begin{gather*}
V'=V-N(S),\\
R'_+ =int ((R_+ -N(S)) \cup S_+ ),\\
\Gamma' =Fr(R'_+ ),\\
R'_- =\partial V' -(R'_+ \cup \Gamma ).
\end{gather*}
Une hi\'erarchie convexe de $(V,\Gamma )$ est une suite de
d\'ecoupages convexes le long de surfaces $\pi_1$--inject\'ees qui
aboutit \`a une union de boules dont les courbes de d\'ecoupage
sont connexes et non vides.

L'ensemble $s(\gamma )$ d'une vari\'et\'e sutur\'ee annulaire
d\'efinit une structure convexe sur un lissage de la vari\'et\'e $V$.
R\'eciproquement, un voisinage tubulaire $\gamma$ des courbes
de d\'ecoupages $\Gamma$ d'une structure convexe d\'efini
une vari\'et\'e sutur\'ee apr\`es introduction de coins le long
de $\partial \gamma$. Une structure convexe est dite {\it tendue}
si la vari\'et\'e sutur\'ee associ\'ee l'est. Il y a de m\^eme une 
correspondance automatique entre hi\'erachies sutur\'ees et 
convexes.

Dans ce contexte, le th\'eor\`eme~\ref{theoreme : hierarchie} a 
une traduction imm\'ediate.

\begin{theoreme}{\rm \cite{HKM1,HKM2}}\qua 
Soit $(V,\Gamma )$ une vari\'et\'e convexe tendue avec\break
$H_2 (V,\partial V;\Z) \neq 0$. Alors $(V,\Gamma )$ admet une
hi\'erarchie convexe tendue {\it bien positionn\'ee}, i.e.\ telle
que toute composante de $\partial S_i$ rencontre $\Gamma$.
\end{theoreme}

\section{Champs de Reeb}

Si $\xi$ est une structure de contact sur une vari\'et\'e
$V$, un champ de Reeb pour $\xi$ est un champ de vecteurs
transversal \`a $\xi$ et dont le flot pr\'eserve $\xi$.

Dans la suite, on suppose que toutes les structures de contact
rencontr\'ees sont coorient\'ees. Tout champ de plans coorient\'e
admet une \'equation globale $\alpha =0$.
On consid\`erera alors toujours des \'equations $\alpha$ de
$\xi$ compatibles avec sa coorientation, c'est-\`a-dire
positives sur tout vecteur direct.
De m\^eme, les champs de Reeb consid\'er\'es
seront toujours suppos\'es positivement transversaux
\`a $\xi$.

Si on fixe une forme de contact $\alpha$ pour $\xi$,
les champs de Reeb $R$ (directs) sont mis en dualit\'e avec
les fonctions $h\co V\rightarrow \R^*_+$ strictement
positives par les \'equations:
$$i_R \alpha =h,\: i_R d\alpha \vert_\xi =-dh\vert_\xi,$$
qui admettent une unique solution $R$.

\`A toute forme de contact $\alpha$ est ainsi associ\'e un champ 
de Reeb (direct pour la coorientation donn\'ee par $\alpha$) 
``naturel", obtenu en prenant $h\equiv 1$.
R\'eciproque\-ment, tout champ de Reeb (direct) pour $\xi$
d\'etermine une unique forme de contact $\alpha$ (positive)
\`a laquelle il est naturellement associ\'e.

La d\'efinition de structure convexe est directement inspir\'ee
de celle de surface convexe, introduite par E~Giroux \cite{Gi1},
et qui d\'ecrit la position d'une surface dans une vari\'et\'e de contact.
Une surface $S\subset (V,\xi )$ est {\it convexe} s'il existe un champ
de vecteurs de contact $X$ transversal \`a $S$. Cette propi\'et\'e 
est g\'en\'erique. Toute surface convexe se d\'ecoupe le long d'une 
multi-courbe $\Gamma_S = \{ x\in S \vert X(x)\in \xi (x)\}$.
La sous-vari\'et\'e $\Gamma_S$ est lisse, transversale \`a $\xi$, et
donc naturellement orient\'ee par la coorientation de $\xi$.
Elle s\'epare $S$ en deux r\'egions $S^+$ et $S^-$, et son orientation
est celle du bord de $S^+$.

Soient $S\subset (V,\xi )$ une surface (\'eventuellement \`a bord)
convexe orient\'ee (et donc coorient\'ee) et $\Gamma$ une courbe 
de d\'ecoupage. On rappelle que $\xi$ est suppos\'ee coorient\'ee.
Un champ de Reeb (direct) $R$ est dit {\it ajust\'e} au
couple $(S,\Gamma )$, s'il est transversal \`a $S - \Gamma$,
positivement sur $S^+$ et n\'egativement sur $S^-$,
et tangent \`a $S$ le long de $\Gamma$, transversal \`a $\Gamma$ 
et rentrant dans $S^+$. Lorsque $R$ est ajust\'e \`a $(S,\Gamma )$,
il existe un voisinage collier $S\times [-\varepsilon ,\varepsilon ]$
de $S\simeq S\times \{0\}$ tel que $R$ soit ajust\'e
\`a $(S \times \{ t\} ,\Gamma \times \{ t\})$ pour tout $t\in 
[-\varepsilon ,\varepsilon ]$.

On convient que si une vari\'et\'e $V$ est orient\'ee,
son bord est orient\'e par la r\`egle: ``la normale sortante
en premier".

Avec cette convention, dans $\R^3$ muni de ses coordonn\'ees
polaires $(r,\theta , z)$ et de la structure de contact
d'\'equation $dz-r^2d\theta =0$, la sph\`ere unit\'e $S$,
orient\'ee comme bord de la boule, est convexe.
Elle est scind\'ee par son \'equateur $\Gamma_S =\{ z=0\}$.
Le champ de Reeb $\frac{\partial}{\partial z}$
pour la forme de contact $dz-r^2 d\theta$ est ajust\'e
au couple $(S,\Gamma_S )$.

Le lemme suivant se d\'emontre en reprenant des arguments
de E~Giroux \cite{Gi1}.             

\begin{lemme} 
Pour toute surface convexe $(S,\Gamma_S ) \subset (V,\xi )$,
il existe un champ de Reeb ajust\'e \`a $\Gamma_S$.
\end{lemme}

Dans la suite, on montrera le th\'eor\`eme suivant:

\begin{theoreme}\label{theoreme : construction}
Si $(V,\Gamma )$ est une vari\'et\'e (connexe) convexe tendue
de bord non vide avec $H_2 (V,\partial V;\Z) \not = 0$, alors 
$(V,\Gamma )$ porte une forme de contact hypertendue dont le 
champ de Reeb est ajust\'e au bord.
\end{theoreme}

Pour rester au plus pr\`es des constructions de Gabai,
on utilise plut\^ot la pr\'esent\-ation ``sutur\'ee", m\^eme
si une \'etude {\it via} les structures convexes est possible.

Soit $(V,\gamma)$ une vari\'et\'e sutur\'ee. Une paire $(\xi, R)$
constitu\'ee d'une structure de contact $\xi$ sur $V$ et d'un 
champ de Reeb $R$ de $\xi$ est {\it adapt\'ee} \`a $(V,\gamma)$ si:

\begin{itemize}

\item $R$ est transversal \`a $R_\pm (\gamma )$, positivement \`a
$R_+ (\gamma )$ et n\'egativement \`a $R_- (\gamma )$.

\item $R$ est tangent \`a $\gamma$.  Les orbites de $R$, de m\^eme
que celles de $\xi \gamma$, feuill\`etent $A(\gamma )=s(\gamma )\times I$
par intervalles $\{ pt\}\times I$, et feuill\`etent chaque composante 
$T$ de $T(\gamma )$ comme la suspension d'un diff\'eomorphisme 
de $S^1$, i.e., le feuilletage donn\'e par $R$ est diff\'eomorphe \`a 
celui donn\'e par les feuilles $\{pt\}\times[0,1]$ de $S^1\times[0,1]$
en identifiant $S^1\times\{0\}$ et $S^1\times\{1\}$ par un diff\'eomorphisme. 
 [Ici, $\zeta \Sigma$ d\'esigne le feuilletage
caract\'eristique induit par la structure de contact $\zeta$ sur 
la surface $\Sigma$.]                                            

\item Toute composante de $\bdry R_+$ est positivement transversale
\`a $\xi$.

\end{itemize}

\s\n
{\bf Mod\`ele dans $\R^3$.}  On consid\`ere $\R^3$ muni
des coordonn\'ees $(r,\theta,z)$ et de la structure
de contact standard $\xi$ donn\'ee par $dz-r^2d\theta=0$.  Soit $(V,\gamma)$
la vari\'et\'e sutur\'ee o\`u $V=D^2\times[-1,1]=\{r\leq 1, -1\leq z\leq 1\}$,
$\gamma=\bdry D^2\times[-1,1]$, et $R_\pm (\gamma)= D^2\times\{\pm 1\}$.
Alors $(\xi,R={\bdry \over \bdry z})$ est adapt\'ee \`a $(V,\gamma)$.

\s
Lorsque $(\xi ,R)$ est adapt\'e \`a la vari\'et\'e $(V,\gamma )$,
et que $\gamma$ est annulaire,
alors $R$ est ajust\'e \`a la structure convexe $(V,\Gamma =s(\gamma ))$
apr\`es un lissage appropri\'e de $\partial V$ le long de $\gamma$.
Ainsi, le th\'eor\`eme~\ref{theoreme : construction} d\'ecoule imm\'ediatement
du r\'esultat suivant:

\begin{theoreme}\label{theoreme : construction2}
Si $(V,\gamma )$ est une vari\'et\'e sutur\'ee tendue de bord non 
vide avec $H_2 (V,\partial V;\Z) \neq 0$, alors $(V,\gamma )$ 
porte une forme de contact hypertendue dont le champ de Reeb 
est adapt\'e au bord.
\end{theoreme}

En particulier, lorsque le bord contient une suture torique $T$,
le champ de Reeb obtenu est tangent \`a $T$. Dans le cas o\`u 
$V$ est bord\'ee par une union non vide de tores, le 
th\'eor\`eme~\ref{theoreme : existence} est donc une cons\'equence
directe de l'application du th\'eor\`eme~\ref{theoreme : construction2}
\`a la vari\'et\'e sutur\'ee $(V,\partial V)$ qui est automatiquement
tendue.

On construit cette forme de contact par collages
successifs le long d'une hi\'erar\-chie sutur\'ee.
\`A chaque \'etape, le champ de Reeb est sans
orbite contractible et adapt\'e aux sutures du bord.

On prend pour convention qu'une structure de contact
sur une vari\'et\'e compacte $V$ est la restriction d'une structure
de contact sur un \'epaississement de $V$.
Le germe de cet \'epaississement est d\'etermin\'e
\`a isotopie relative au bord pr\`es par sa trace sur le bord, ce qui enl\`eve
toute ambigu\"\i t\'e sur l'\'epaississement consid\'er\'e.
Dans les paragraphes qui suivent, on peut \^etre amen\'e
\`a consid\'erer des d\'eformations de $\partial V$,
ce qui signifie des d\'eformations dans ce germe d'\'epaississe\-ment.

\section{Lemme de flexibilit\'e}

Le lemme suivant est essentiellement contenu
dans \cite{Gi1}. C'est la pierre angulaire de notre construction.

\begin{lemme}\label{lemme : flexibilite}{\rm\cite{Gi1}}\qua
Soit $S$ une surface compacte avec $\bdry S\not=\emptyset$, et $dt+\beta$
une 1--forme de contact sur $S\times \R$, o\`u $S\times\R$ a pour coordonn\'ees $(x,t)$ et $\beta$
est le rappel d'une 1--forme de $S$.  Si $\beta'$ est la primitive d'une forme d'aire $d\beta'$
sur $S$ et $\beta|_{\bdry S}=\beta'|_{\bdry S}$, alors il existe un diff\'eomorphisme
$\phi=(\phi_1,\phi_2)\co  S\times \R\rightarrow S\times \R$, o\`u $\phi_1(x,t)
=\phi_1(x,t')$, $\phi_2(x,t+t')=\phi_2(x,t)+t'$, $\phi^*(dt+\beta')=dt+\beta$, et $\phi=id$
sur $\bdry S\times \R$.
\end{lemme}

\begin{proof}
[D\'emonstration]
On suit les lignes d'un argument de Moser.  Pour $s\in[0,1]$, soit
$\alpha_s=dt+\beta_s$, o\`u $\beta_s$ est l'interpolation $(1-s)\beta+s\beta'$.
On observe que $d\beta_s$ est une forme d'aire pour tout $s$, et donc que $\alpha_s$ est une
forme de contact pour tout $s$.  On cherche \`a r\'esoudre l'\'equation d'inconnue $X_s$:
\begin{equation}  \label{Lie}
\mathcal{L}_{X_s} \alpha_s={d\alpha_s\over ds}=\beta'-\beta.
\end{equation}                       
On \'ecrit $X_s=f_s {\bdry\over \bdry t} +Y_s$ avec $f_s$ une fonction sur $S$ et
$Y_s$ un champ de vecteurs sur $S$, i.e., $f_s$ et $Y_s$ sont ind\'ependants de $t$.
Si on d\'eveloppe par la formule de Cartan, il suffit de r\'esoudre:
\begin{equation}
i_{Y_s} d\beta_s=\beta'-\beta, \alpha_s(X_s)=f_s+\beta_s(Y_s)=0.
\end{equation}
On obtient $Y_s$ \`a l'aide de la premi\`ere \'equation, et $f_s$ avec la seconde.
En int\'egrant $X_s$, on obtient le diff\'eomorphisme recherch\'e $\phi$.  De plus,
comme $\beta=\beta'$ sur $\bdry S$, $Y_s=0$ et $f_s=0$ sur $\bdry S$, et
donc $\phi=id$ sur $\bdry S\times\R$.
\end{proof}

On peut utiliser ce lemme de flexibilit\'e dans la situation suivante:

Soit $(V,\gamma)$ une vari\'et\'e sutur\'ee avec une paire adapt\'ee $(\xi,R)$, et $S$
une composante connexe de $R_+(\gamma)$.  (L'argument pour $R_-(\gamma)$ est
similaire.)  Un voisinage de $S$ peut \^etre plong\'e dans $S\times \R$ de sorte que
$S$ s'envoie sur $S\times\{0\}$, $R$ s'envoie sur ${\bdry\over \bdry t}$, et
la 1--forme de contact correspondante s'envoie sur
$dt+\beta$.  Si on veut remplacer $\beta$ par $\beta'$ qui co\"\i ncide avec
$\beta$ le long de $\bdry S$ et telle que $d\beta'$ soit une forme d'aire, on
applique le lemme de flexibilit\'e  pour obtenir
$S'=\phi^{-1}(S\times\{t\})\subset (S\times\R, dt+\beta)$.  Si $t$
est assez grand, alors $S'$ et $S=S\times\{0\}$ sont disjoints et
cobordent une r\'egion $V'$ diff\'eomorphe \`a $S\times I$.  
L'attachement de $V'$ \`a $V$ le long de $S$ a pour effet de modifier 
le feuilletage caract\'eristique de $S$ de celui donn\'e par $\beta$ 
\`a celui donn\'e par $\beta'$.

\section{Collage}

Soit $(V,\gamma)\stackrel {S}\rightsquigarrow (V',\gamma')$ une
d\'ecomposition de vari\'et\'e sutur\'ee annulaire tendue bien
positionn\'ee, i.e.\ on suppose que $\bdry S$ est non vide, que $S$
est connexe et $\pi_1$--inject\'ee dans $V$, que $(V,\gamma)$ a des
sutures annulaires et que chaque composante de $\partial S$ rencontre 
$\gamma$ le long d'arcs. On suppose \`a pr\'esent que $(\xi',R')$ est 
adapt\'e \`a $(V',\gamma')$.

On construit la paire $(\xi,R)$ adapt\'ee \`a $(V,\gamma )$ en collant 
$S'_+\subset R'_+(\gamma')$ \`a $S'_-\subset R'_-(\gamma')$, o\`u 
$S'_\pm$ sont les copies positives et n\'egatives de $S$ obtenues 
par d\'ecoupage le long de la surface orient\'ee $S$. Comme toutes 
les composantes de $\bdry S$ intersectent non trivialement $\gamma$, 
toute composante du bord orient\'e $\bdry S'_+$ est une courbe 
ferm\'ee lisse par morceaux qui est constitu\'ee d'une union d'arcs 
lisses orient\'es $a_1^+,b_1^+,\dots, a_k^+,b_k^+$
o\`u (i) la fin de $a^+_i$ est le d\'ebut de $b^+_i$ et la fin de $b^+_i$
est le d\'ebut de $a^+_{i+1}$ ($i+1$ est pris  modulo $k$), (ii) $\gamma'\cap \bdry S'_+=\cup_{i=1}^k a^+_i$, et
(iii) $k\geq 1$.  De fa\c con similaire, toute composante de  $\bdry S'_-$ est compos\'ee
d'arcs $a^-_1, b^-_1,\dots,a^-_k,b^-_k$, o\`u (i)  la fin de  $a^-_i$ est le d\'ebut de
$b^-_i$ et la fin de $b^-_i$
est le d\'ebut de $a^-_{i+1}$ (mod $k$), (ii) $\gamma'\cap \bdry S'_-=\cup_{i=1}^k b^-_i$, et
(iii) $k\geq 1$.

En prenant un diff\'eomorphisme d'un voisinage de $R'_+(\gamma')$ dans $R'_+(\gamma')
\times \R$ envoyant $R'_+(\gamma')$ sur $R'_+(\gamma')\times\{0\}$,
on peut supposer que $R'={\bdry\over
\bdry t}$ et que la forme de contact est $dt+\beta_+$.

\s\n
{\bf \'Etape 1}

\begin{affirmation}
Il existe une 1--forme $\beta'_+$ sur $R'_+(\gamma')$ telle que $d\beta'_+$
soit une forme d'aire sur $R'_+(\gamma')$,
que $\beta'_+=\beta_+$ sur $\bdry R'_+(\gamma')$, et que
$\bdry S'_+$ soit une courbe transversale positive pour $dt+\beta'_+$.
De plus, $\beta'_+(a^+_i)
=\beta'_+(b^+_i) =\varepsilon$ pour un petit $\varepsilon>0$.   (Ici, on peut
\^etre amen\'e \`a r\'etr\'ecir $S'_+$
si n\'ec\'essaire.)
\end{affirmation}

\begin{figure}[ht!]
        {\epsfxsize=4in\centerline{\epsfbox{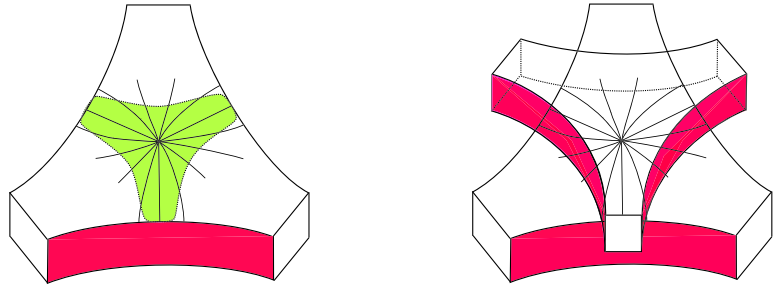}}}
        \caption{Pr\'eparation du feuilletage caract\'eristique
        de  $S'_+$ et collage de $S'_+$ et $S'_-$}
        \label{collage}
\end{figure}

\begin{proof}
[D\'emonstration]
On consid\`ere un squelette legendrien $K$ pour $S'_+$ de la fa\c con suivante.
On d\'efinit $\beta'_+$ au voisinage d'un point $p$ dans l'int\'erieur
de $S'_+$ de
sorte que $p$ devienne une singularit\'e elliptique positive (par exemple,
$\beta'_+={1\over 2}(xdy-ydx)$ pr\`es de $(0,0)$), et on attache successivement
des anses d'indice $1$, i.e., on prend des singularit\'es hyperboliques positives
(donn\'ees par exemple par $\beta'_+= 2xdy + ydx$ pr\`es de $(0,0)$), et on connecte
leurs s\'eparatrices instables \`a $p$. L'union $K'$ du point elliptique,
des points hyperboliques et de leurs s\'eparatrices instables
aura le m\^eme type topologique que
$S'_+$.  Finalement, on ajoute un arc legendrien de $p$ \`a $a^+_i$, $i=1,\dots,k$, pour
obtenir $K$.  Un petit voisinage convenable de $K$ (de bord transversal
au feuilletage caract\'eristique dirig\'e par $\ker \beta'_+$) sera
$S'_+$. Comme $\beta'_+$ est nulle le long de $K$,
on peut imposer n'importe qu'elle $\beta'_+$--longueur assez
petite aux arcs $a^+_i$ et $b^+_i$.

On note \'egalement que la $\beta'_+$--longueur d'une courbe 
ferm\'ee $l$ (ou d'un arc) transversale \`a la structure $\xi$ d\'etermine 
$\beta'_+$ \`a diff\'eomorphisme pr\`es au voisinage de $l$ (relativement 
au bord de l'arc), comme l'explicite le lemme suivant.

\begin{lemme}
Soit $S$ une surface \`a bord et $dt+\beta_0$ et $dt+\beta_1$
deux formes de contact d\'efinies au voisinage $S\times [-1,1]$ 
de $S=S\times \{ 0\}$, o\`u $t$ est la coordonn\'ee dans $[-1,1]$ 
et $\beta_0$ et $\beta_1$ sont les rappels de 1--formes sur $S$. 
On suppose que $l \subset \partial S$
est une courbe ferm\'ee ou un arc transversal \`a $\ker \beta_i$, $i=0,1,$ et que 
$\int_l \beta_0 =\int_l \beta_1$. Il existe un diff\'eomorphisme de 
$S$ qui conjugue $\beta_0$ et $\beta_1$ pr\`es de $l$.
\end{lemme} 

\begin{proof}[D\'emonstration]    
On suppose que $l$ est une courbe ferm\'ee.
Un voisinage de $l\subset \partial S$ est conjugu\'e \`a $\R /\Z 
\times [0,1]=\{ (x,y) \}$, $l=\R /\Z \times \{ 0\}$, avec $\beta_1 
=f(x,y)dx$. La fonction $f$ est strictement positive et
l'application $(x,y)\mapsto (\int_0^x f(u,0)du ,y)$ 
est un diff\'eomorphisme de $\R /\Z \times [0,1]$
dans $\R /L_l \Z \times [0,1]$ qui envoie $\beta_1$ sur
$F(x,y)dx$, $F(x,0)=1$. On pose $L_l= \int_0^1 f(x,0)dx$.
La condition de contact dit que ${\partial F\over \bdry_y} >0$.
Le diff\'eomorphisme $(x,y)\mapsto (x,F(x,y))=(X,Y)\in \R /L_l \Z
\times [1,1+\varepsilon]$
envoie $F(x,y)dx$ sur $YdX$.
On peut donc conjuguer les formes $\beta_0$ et $\beta_{1}$ au m\^eme
mod\`ele, et donc entre elles.
\end{proof}

Il reste \`a \'etendre $\beta'_+$ de  $S'_+$ \`a $R'_+(\gamma')$, 
en utilisant le fait que les $\beta'_+$--longueurs des $b^+_i$ 
peuvent \^etre choisies arbitrairement petites. On d\'ecrit 
succinctement cette op\'eration.

On note $R_0$ une composante connexe de $R'_+ (\gamma' ) 
\setminus int(S'_+ )$. La surface $R_0$ poss\`ede des coins dans 
son bord. Comme dans \cite{Gi1}, on peut munir $R_0$ d'une 
$1$--forme ~$\beta_0$ avec les propri\'et\'es suivantes:
\begin{itemize}
\item $d\beta_0$ est une forme de surface sur $R_0$;
\item $\ker \beta_0$ est n\'egativement transversal \`a
$\partial R_0 \cap (\cup_i b^+_i )$ et positivement
trans\-versal \`a $\partial R_0 \setminus (\cup_i b^+_i )$.
\item la $\beta_0$--longueur de chaque composante de
$\partial R_0 \setminus (\cup_i b^+_i )$ est inf\'erieure \`a
sa $\beta'_+$--longueur.
\end{itemize}
On choisit alors $\beta'_+$ de sorte que la $\beta'_+$--longueur
de chaque composante de $\partial R_0 \cap (\cup_i b^+_i )$
soit inf\'erieure \`a sa $\beta_0$--longueur.

Dans cette situation, on peut modifier (comme, par exemple, dans le lemme~\ref{lemme : longueur})
 $\beta_0$ pr\`es de chaque composante  de  $\partial R_0 \cap (\cup_i b^+_i )$
et de $\partial R_0 \setminus (\cup_i b^+_i )$ pour faire co\"\i ncider
leurs $\beta_0$-- et $\beta'_+$--longueurs.
La forme $\beta'_0$ obtenue en d\'eformant $\beta_0$ est l'extension
de $\beta'_+$ recherch\'ee.

Il faut encore s'occuper des coins de $S'_+$. Ce sera fait dans 
la prochaine \'etape.   Le diagramme de gauche dans la figure~\ref{collage}
d\'ecrit la fa\c con dont le feuilletage caract\'eristique de $S'_+$
est pr\'epar\'e avant le collage.
\end{proof}                            

\s\n
{\bf \'Etape 2}\qua (Forme normale pr\`es des coins)\qua
Pr\`es de la fin $q$ de $b^+_i$ (qui est le d\'ebut de $a^+_{i+1}$),
on normalise $\beta'_+$ comme suit: on plonge un voisinage $N(q)$ de $q$ dans
$R'_+(\gamma')$ sur $\R^2$ avec ses coordonn\'ees $(x,y)$, en envoyant $q$ sur $(0,0)$,
$\bdry R'_+(\gamma') \cap N(q)$ sur la droite $y=x$, et
$\beta'_+$ sur $-(y+1) dx$. En modifiant $S'_+$ si n\'ec\'essaire, on peut supposer que
$b^+_i$ est envoy\'e sur $y=-x$ et que $a^+_{i+1}$ est envoy\'e sur $y=x$.
De fa\c con similaire, pr\`es de la fin de  $a^+_i$, on envoie $\bdry R'_+(\gamma')\cap N(q)$
sur la droite $y=-x$, $a^+_i$ sur $y=-x$, et $b^+_i$ sur $y=x$.

\s\n
{\bf \'Etape 3}\qua (Lissage)\qua Gr\^ace \`a la deuxi\`eme \'etape,
on peut s'arranger pour que $\beta'_+$ sur $S'_+$ et
$\beta'_-$ sur $S'_-$ co\"\i ncident pr\`es des coins. Maintenant, 
par la premi\`ere \'etape, on a rendu les $\beta'_\pm$--longueurs des 
$a^\pm_i$ et $b^\pm_i$ toutes \'egales \`a $\varepsilon$, si
bien qu'on peut faire co\"\i ncider $\beta'_\pm$ sur $\partial S'_\pm$ 
relativement \`a leurs coins.  On utilise \`a pr\'esent le lemme de flexibilit\'e
pour obtenir une tranche diff\'eomorphe \`a
$R'_\pm(\gamma')\times[0,1]$ que l'on ajoute \`a $(V',\gamma')$ le long de $R'_\pm(\gamma')$.
On peut alors supposer qu'il existe un diff\'eomorphisme entre $S'_+$ et $S'_-$
qui envoie $a^+_i$ sur $a^-_i$ et $b^+_i$ sur $b^-_i$ et qui \'echange les formes de contact.

On \'epaissit $(V',\gamma')$ en lui attachant $S\times [0,1]$ avec une structure de contact
$[0,1]$--invariante et un champ de Reeb ${\bdry \over \bdry t}$.  On prend alors
$V'\cup (S\times[0,1])$, o\`u $S\times\{0\}$ est attach\'e \`a $S'_+$ et
$S\times\{1\}$ est attach\'e \`a $S'_-$.  Soit maintenant
$P_\pm=\overline{R'_\pm(\gamma')\setminus S'_\pm}$.
On prend $P_\pm\times[0,1]$ avec une structure de contact $[0,1]$--invariante
et un champ de Reeb ${\bdry \over \bdry t}$.   On colle $P_+\times\{0\}$ \`a $P_+$,
$P_-\times\{1\}$ \`a
$P_-$, et $\bdry P_\pm\times[0,1]$ \`a $(\bdry S\cap P_\pm) \times [0,1]$ de $S\times[0,1]$.
On note $V_1$ la vari\'et\'e obtenue.

Il reste \`a arrondir les coins.   Soit $A=S^1\times[0,1]$ la composante connexe de $A(\gamma')$
avec $a_i^+\subset S^1\times\{1\}$, et $A\times [0,\varepsilon [$ son voisinage dans $V'$, o\`u
$A\times \{0\}=A$. Les coordonn\'ees sur $S^1\times[0,1]\times[0,\varepsilon [$
seront not\'ees $(\theta,x,y)$,
et on prendra $a_i^+=\{ x=1, y=0, -{\pi\over 2} \leq \theta\leq {\pi\over 2}  \}$.
Soit $f\co [0,\varepsilon [\rightarrow [0,1]$ une fonction lisse
d\'ecroissante avec $f(y)=1$  pr\`es de $y=0$ et
$f(y)=0$ pr\`es de $y=\varepsilon$. On soustrait \`a $V'$ les ensembles
$\{x\leq f(y), -{\pi\over 2} \leq \theta\leq {\pi\over 2} \}$, $\{ x\leq f(\sqrt {y^2+(\theta\pm{\pi\over 2})^2})\}$.

On traite \`a pr\'esent le cas des sutures toriques. Il existe alors des 
paires de composantes de bord $\delta_\pm \subset\bdry S'_\pm$ 
qui cobordent une composante connexe $A_{\delta_\pm}$ de  
$\gamma'$ et qui sont identifi\'ees par l'application de collage pour 
donner une suture torique. Comme d'habitude, les orientations
de $\delta_\pm$
 sont induites par celles de 
$S'_\pm$. Dans ce cas, on d\'efinit $S'_+$ comme
\'etant le voisinage de $K$ comme dans la premi\`ere \'etape, \'etendu 
des composantes de $S'_+\setminus K$ dont une des composantes 
de bord est $\delta_+$. Comme le champ $R'$, associ\'e \`a la forme 
de contact $\alpha'$ d\'efinie sur $V'$, est tangent \`a $A_{\delta_\pm}$,
on a (par la formule de Stokes) $\int_{\delta_-} \alpha' =\int_{\delta_+} 
\alpha' >0$. Pour pouvoir appliquer la construction pr\'ec\'edente,
il faut encore \^etre s\^ur que les $\alpha'$--longueurs des diff\'erentes 
paires $\delta_\pm$ sont \'egales entre elles.  Pour cela, on constate 
qu'on peut toutes leur faire prendre une valeur assez grande fix\'ee 
en \'epaississant $V'$ le long des sutures $A_{\delta_\pm}$:
on colle \`a $A_{\delta_\pm}$ un produit $\R /\Z \times [-1,1]
\times [a,b]$, muni de coordonn\'ees $(x,y,t)$, en identifiant
$A_{\delta_\pm}$ \`a $\{ t=a\}$ \`a l'aide du lemme suivant:

\begin{lemme}[Ajustement des longueurs]  \label{lemme : longueur}
Soit $\xi =\ker (\cos t dx -\sin t dy)$ sur $W=\R /\Z \times [-1,1]
\times [a,b]$, muni de coordonn\'ees $(x,y,t)$, o\`u
${\pi\over 2}< a<b< \pi$. Soit de plus $R$ un champ
de Reeb pour $\xi$, d\'efini
pr\`es de $\{ t=a\}$, tangent \`a $\{ t=a\}$ et transversal \`a $\R /\Z
\times \{ y\} \times [a,b]$ pour tout $y\in [-1,1]$. Il existe $L_0 >0$
tel que pour tout $L>L_0$, il existe une extension de $R$ en un
champ de Reeb, not\'e \'egalement $R$, sur $W$ tel que: $R$
soit transversal \`a chaque $\R /\Z \times \{ y\} \times [a,b]$, $R$ soit
tangent \`a $\{ t=b\}$, et $\int_{\R /\Z \times \{ \pm 1\} \times \{b\}}
\alpha_R =L$.  
\end{lemme}

\begin{proof}[D\'emonstration] 
Le champ $R=R_x{\bdry \over \bdry x} + R_y {\bdry \over \bdry y}
+R_t{\bdry \over \bdry t}$ est d\'etermin\'e par la donn\'ee  d'une
fonction g\'en\'eratrice $H\co  W\rightarrow \R^*_+$ qui v\'erifie:
\begin{gather*}
R_x = H \cos t -{\bdry H\over \bdry t}\sin t\\ 
R_y = - H \sin t -{\bdry H\over \bdry t} \cos t\\
R_t = -{\bdry H\over \bdry x} \sin t  -{\bdry H\over \bdry y}\cos t.
\end{gather*}
La fonction $H$ est donn\'ee pr\`es de $\{ t=a\}$. On l'\'etend
en une fonction not\'ee \`a nouveau $H$ avec ${\bdry H\over \bdry t} 
<0$ en dehors d'un petit voisinage de $\{ t=a\}$ (ce qui assure que 
$R_y <0$) et telle que ${\bdry H\over \bdry x} ={\bdry H \over \bdry y} 
=0$ le long de $\{ t=b\}$ (ce qui assure que $R$ est tangent
\`a $\{ t=b\}$).

La forme de contact associ\'ee \`a un tel $R$ est ${1\over H} 
(\cos tdx-\sin t dy)$. La longueur $L$ du bord $\{ y=\pm 1, t=b\}$ est
$ {-\cos b\over H|_{\{t=b\}}}$. Ici, $H$ est constant sur
$\{t=b\}$. On peut choisir $H\vert_{\{ t=b\}}$ arbitrairement petit.
\end{proof}

\begin{lemme}\label{lemme : contractible} 
Si  $R'$ n'a pas d'orbite p\'eriodique contractible dans $V'$, 
alors $R$ n'a pas d'orbite p\'eriodique contractible dans $V$.
\end{lemme}

\begin{proof}[D\'emonstration] 
Soit $\mathcal{P}$ une orbite p\'eriodique de $R$. Le champ de 
Reeb $R$ est positivement transversal \`a la surface de d\'ecoupage 
orient\'ee $S$. Si $\mathcal{P}$ rencontre $S$, alors, l'intersection 
de $\mathcal{P}$ et de $S$ en homologie est non nulle, \'egale au 
nombre de points d'intersection de $\mathcal{P}$ et de $S$, et donc 
$\mathcal{P}$ est non contractible dans $V$. Si $\mathcal{P}$ ne 
rencontre pas $S$, alors $\mathcal{P}$ est contenue dans $V'$.
Comme la surface de collage $S$ est $\pi_1$--inject\'ee dans $V$,
le groupe fondamental de $V'$ s'injecte dans celui de $V$
d'apr\`es le th\'eor\`eme de Seifert--Van Kampen.
L'orbite $\mathcal{P}$ qui n'\'etait pas contractible dans $V'$
ne l'est toujours pas dans $V$.
\end{proof}

\section{Preuve des th\'eor\`emes}

Sous les hypoth\`eses du th\'eor\`eme~\ref{theoreme :
construction2}, la vari\'et\'e $(V,\gamma )$ poss\`ede une
hi\'erarchie sutur\'ee qui est annulaire, tendue, et bien
positionn\'ee apr\`es un premier d\'ecoupage (le long d'une
surface $S$ \'even\-tuellement non connexe et qui intersecte
$A(\gamma )$ le long d'arcs) \'eliminant les sutures toriques
(voir \cite{HKM1} et le th\'eor\`eme~\ref{theoreme : hierarchie}).
En partant du mod\`ele de champ de Reeb d\'ecrit sur la boule
sutur\'ee $D^2 \times [0,1]\subset \R^3$ dans le paragraphe~4
(cf.\ {\it Mod\`ele dans $\R^3$}), on
obtient par collages successifs une paire $(\xi ,R)$ adapt\'ee \`a
$(V,\gamma )$. D'apr\`es le lemme~\ref{lemme : contractible},
appliqu\'e \`a chaque \'etape du collage, le champ $R$ est sans
orbite p\'eriodique contractible.

Le feuilletage tendu construit par Gabai sur $V$ est ``orthogonal"
au champ de Reeb $R$. On part du feuilletage horizontal en disque 
de $D^2 \times [0,1]$, alors que le champ de Reeb \`a cette \'etape 
est tangent \`a la direction verticale. \`A chaque \'etape, on colle \`a 
la vari\'et\'e $V_i$ des produits, dont la direction verticale est 
toujours donn\'ee par le champ de Reeb, alors que le feuilletage
s'\'etend transversalement \`a cette verticale. Le fait que le champ 
de Reeb $R$ soit transversal au feuilletage tendu obtenu par Gabai 
dans \cite{Ga1} est donc automatique.

\begin{remarque}
{\rm L'argument essentiel est que le champ de Reeb est transversal 
\`a chaque surface de recollement. Il n'est pas possible d'obtenir 
cette propri\'et\'e lorsque la surface de recollement est sans bord:
la forme $d\alpha$ serait alors une forme de surface exacte, contredisant
le th\'eor\`eme de Stokes. C'est pourquoi cette construction ne marche
en g\'en\'eral que sur les vari\'et\'es de dimension trois \`a bord.}
\end{remarque}

\begin{remarque}\label{remarque : voisinage}
{\rm Dans le th\'eor\`eme~\ref{theoreme : existence}, on peut obtenir 
en fait mieux qu'un champ de Reeb tangent au bord: un voisinage 
du bord est feuillet\'e par des tores satur\'es par le flot.}
\end{remarque}

\section{Contr\^ole de la pente au bord}

On prouve ici un corollaire du th\'eor\`eme~\ref{theoreme : existence}
qui permet, la direction du champ de Reeb \'etant fixe, de contr\^oler 
la pente du feuilletage caract\'eristique du bord.

Soit $T$ un tore de dimension $2$ orient\'e et $[m]$, $[l]$
$\in H_1 (T; \Z )$ une base de $H_1 (T;\Z )$, avec $i([m] ,[l])=1$.
Les cycles asymptotiques de Schwartzman \cite{Sc} permettent
d'associer \`a tout feuilletage non singulier $\F$ de $T$ une
direction dans $H_1 (T;\R )=\R [m] +\R [l]$, dont le coefficient
directeur est appel\'e la {\it pente} de $\F$. Par exemple, la
pente d'un feuilletage de $T$ par des cercles homologues \`a $[m]$
vaut $0$.

Soit $V$ une vari\'et\'e compacte, irr\'eductible et orient\'ee
bord\'ee par une union non vide de tores $T_1$,...,$T_n$. On
suppose que $S$ est une surface minimale (i.e.\ sans composantes
closes et qui minimise la norme de Thurston dans $H_2(V,\bdry V;\Z)$)
orient\'ee dans $V$ qui
rencontre tous les $T_i$ le long d'une famille de courbes
homologues dans $T_i$. Chaque $T_i$ est orient\'e comme bord de
$V$. On fixe alors, pour $1\leq i\leq n$ une base $[m_i]$, $[l_i]$
de $H_1 (T_i ;\Z )$ comme indiqu\'e ci-dessus, avec de surcro\^\i
t, $[m_i]=[s_i]$, o\`u $s_i$ est une composante de $S\cap T_i$. La
pente d'un feuilletage de $T_i$ est alors calcul\'e dans la base
$([m_i ],[l_i ])$.

\begin{corollaire}\label{corollaire : controle}
Soit $V$ une vari\'et\'e compacte, orientable, irr\'eductible,
bord\'ee par une union non vide de tores $T_1$,...,$T_n$. On note
$S$ une surface minimale dans $V$ qui rencontre toutes les
composantes de $\partial V$. Pour tout $\varepsilon <0$
, 
il existe une forme de contact
hypertendue $\alpha_{\varepsilon}$ sur $V$, dont le champ de Reeb est
tangent \`a $\partial V$ et transversal \`a $S$, de direction
ind\'ependante de $\varepsilon$, et tel que la pente de chaque
feuilletage $\xi_\varepsilon \partial T_i$ ($\xi_\varepsilon =\ker
\alpha_\varepsilon$) soit dans l'intervalle
$[\varepsilon,0 [$.
\end{corollaire}

\begin{proof}[D\'emonstration]
Soit $\alpha$ une forme de contact adapt\'ee \`a $(V\setminus S,
\gamma= (\cup_i T_i)\setminus \bdry S)$, donn\'ee par le th\'eor\`eme
\ref{theoreme : construction2} et $R$ son champ de Reeb. 
On se donne alors un petit \'epaississement $U=S\times [-a,a]$ de $S$
avec $R={\bdry\over \bdry t}$, o\`u $t$ est la variable de $[-a,a]$.

On pose $\alpha_\varepsilon =\alpha$ en dehors de $U$.
Sur $U$, on remplace $\alpha$ par la forme
de contact $\alpha_\varepsilon =\alpha -{1\over \varepsilon} \chi (t)dt$,
o\`u $\chi \co [-a,a] \rightarrow \R$ vaut $0$ en $\pm a$
et est strictement positive \`a l'int\'erieur.
Quand $\varepsilon$ tend vers $0$, la structure $\ker \alpha_\varepsilon$
tend vers le champ de plan horizontal $TS$, et
donc la pente du feuilletage trac\'e sur
chaque composante de bord tend vers $0$.
En tout point de $V$, $d\alpha_\varepsilon =d\alpha$ et donc le
champ de Reeb associ\'e \`a $\alpha_\varepsilon$ est proportionnel 
\`a celui donn\'e par $\alpha$.
\end{proof}

\section{Preuve du corollaire \ref{corollaire : toroidal} --
Obturations de Dehn}

Il s'agit ici de d\'emontrer que la forme de contact
produite par les th\'eor\`emes \ref{theoreme : existence}
et \ref{corollaire : controle}
s'\'etend de mani\`ere contr\^ol\'ee \`a des
recollements de tores \'epaissis $(T^2 \times [a,b] ,\cos tdx-\sin tdy=0 )$
ou des obturation de Dehn.

Le feuilletage trac\'e sur chaque tore du bord
est laiss\'e invariant par le flot
du champ de Reeb. Pour $T\subset \partial V$, le feuilletage
$\xi T$ est lin\'eaire:
soit il est de pente rationnelle et il poss\`ede
une orbite p\'eriodique, n\'ec\'essairement
transversale au champ de Reeb, et il s'agit alors d'un feuilletage
en cercles, sinon c'est un feuilletage de pente irrationnel,
l\`a encore lin\'eaire.

\subsection{Preuve du corollaire~\ref{corollaire : toroidal}}

On consid\`ere la structure de contact $\xi$ d\'efinie
sur $U=T^2 \times [a,a' ]$,
muni des coordonn\'ees $((x,y),t)$,
par l'\'equation $\alpha =\cos tdx -\sin tdy$.
Soit $p \co  T^2 \times [a,a' ] \rightarrow T^2$
la projection de fibre $[a,a' ]$.
On note $R_0$ le champ de Reeb associ\'e \`a $\alpha$.

\begin{lemme}\label{lemme : obturation1}
Soit $R_1$ un germe de champ de Reeb pour $\xi$
pr\`es du tore $\{ t=a \}$, dont la composante
sur le vecteur ${\bdry\over \bdry t}$ est nulle pr\`es de $\{ t=a\}$.
Il existe $c >0$, et une extension $R$ de $R_1$ sur $U$
qui vaut $R_0$ sur $T^2 \times [a+c,a']$ et qui est partout
(positivement) transversale \`a un conjugu\'e,
relativement \`a $T^2 \times \{ a,a+c\}$, du feuilletage de codimension
un d'\'equation $\cos{a} dx -\sin{a} dy=0$ sur $T^2 \times [a,a+c]$.
\end{lemme}                   

\begin{proof}[D\'emonstration]  
On fait d'abord la remarque suivante.
Un champ de vecteurs $R$ sur $U$ est un champ de Reeb
si et seulement s'il existe une fonction
$H \co U \rightarrow \R^*_+$ telle que:
\begin{gather*}
R_x = H \cos t -{\bdry H\over \bdry t}\sin t \\
R_y = - H \sin t -{\bdry H\over \bdry t} \cos t \\
R_t = -{\bdry H\over \bdry x} \sin t  -{\bdry H\over \bdry y}\cos t.
\end{gather*}
La troisi\`eme \'equation implique que
si $R_t \equiv 0$, la fonction $H$ est constante
sur les feuilles du feuilletage caract\'eristique $\xi (T^2 \times \{ t\})$.
Lorsque la pente de ce  feuilletage est irrationnelle, ses feuilles sont denses
et donc $H\vert_{T^2 \times \{ t\}}$ est constante.
Cette situation se produit pour un ensemble dense de valeurs
du param\`etre $t$, et donc ${\bdry H\over \bdry x} =
{\bdry H\over \bdry y} =0$.

Dans le contexte du lemme~\ref{lemme : obturation1},
on en d\'eduit que $R_1$ ne d\'epend que de $t$.

Il existe $\pi >c>0$ tel que la base $(p_* R_0 (\{ t=a+c\} )=
\cos (a+c){\bdry \over \bdry x} -\sin (a+c){\bdry\over \bdry y}$, 
$p_* R_1 (\{ t=a\} ))$
soit directe et que $R_0 (\{ t=a+c\} )$ soit positivement
transversal au feuilletage $\{ \cos a  dx -\sin a dy=0\}$.

Comme $R_1$ est tangent \`a $\{ t=a \}$,
l'\'equation associ\'ee \`a $R_1$ en $\{ t=a\}$
s'\'ecrit: $udx -vdy =0$.
On peut alors trouver $f,g \co [a,a+c] \rightarrow \R$
telles que:
\begin{itemize}
\item $f(a)=u$, $g(a)=v$, $f(a+c)=\cos (a+c)$, $g(a+c)=\sin (a+c)$;
\item $\forall t\in [a,a+c]$, $(f(t),g(t))\neq (0,0)$;
\item $\forall t\in [a,a+c]$, $f(t)g'(t)-f'(t)g(t)>0$ (condition de contact
pour $f(t)dx-g(t)dy$);
\item le vecteur de coordonn\'ees $(f(t),g(t))$ soit
positivement transversal au feuilletage de $\R^2$ d'\'equation
$\cos a dx-\sin a dy=0$.
\end{itemize}

Sous ces conditions, la forme d\'efinie sur $T^2 \times [a,a+c]$ comme
$f(t) dx-g(t)dy$ est de contact. Son champ de Reeb est $f(t)
{\bdry\over \bdry x}
-g(t){\bdry\over\bdry y}$. Elle est l'interpolation recherch\'ee.
Sur $T^2 \times [a,a+c]$, on a remplac\'e $\xi$
par la structure $\ker (fdx-gdy)$ qui lui est conjugu\'ee
relativement au bord.
\end{proof}

Pour montrer le corollaire~\ref{corollaire : toroidal},
on se donne un tore $T\subset V$ incompressible et
on d\'ecoupe $V$ le long de $T$ pour obtenir une
vari\'et\'e \`a bord $V'$.
Le th\'eor\`eme~\ref{theoreme : existence} permet de munir
$V'$ d'une forme de contact hypertendue $\alpha'$
dont le champ de Reeb $R'$ est tangent \`a $\partial V'$ et
transversal \`a la premi\`ere surface de d\'ecoupage de
$V'$. On a m\^eme (remarque~\ref{remarque : voisinage})
un voisinage du bord feuillet\'e par des tores satur\'es
par le champ $R'$.
On colle alors \`a $V'$ un tore \'epais
$(T^2 \times [a,b] ,\cos tdx -\sin t  dy=0 )$ pour obtenir une
vari\'et\'e de contact diff\'eomorphe \`a $V$.
On \'etend la forme $\alpha'$ pr\`es de $a$ et $b$ \`a l'aide
du lemme~\ref{lemme : obturation1} pour qu'elle donne
la m\^eme \'equation $\cos tdx -\sin tdy$
en $T^2 \times \{ a+c \}$ et $T^2 \times \{b -c' \}$.

On \'etend alors la forme obtenue par $\cos tdx-\sin tdy$
sur $T^2 \times [a+c ,b-c' ]$.

On obtient finalement ainsi une version contr\^ol\'ee
du corollaire~\ref{corollaire : toroidal}, 
par combinaison de la preuve pr\'ec\'edente et de celle
du corollaire~\ref{corollaire : controle}: comme $\partial V'$
est feuillet\'e par des orbites de $R$, toute orbite
p\'eriodique  reste soit dans $V'$, soit dans son compl\'ementaire.

Soit $V$ une vari\'et\'e toro\"\i dale irr\'eductible, orientable
et close et $T\subset V$ un
tore incompressible. On note $U =T^2 \times [a,b]$ un voisinage
tubulaire de $T$ et $V' =V\setminus T^2 \times ]a,b[$.
On note de plus $S\subset V'$ la premi\`ere
surface de d\'ecoupage de $V'$, qui rencontre les deux composantes
de $\partial V'$.

\begin{theoreme}\label{theoreme : tore}
Pour tout $\varepsilon <0$, il existe une forme de contact
hypertendue $\alpha_\varepsilon$ sur $V$ et deux r\'eels $c,c'>0$
avec les propri\'et\'es suivantes:
\begin{itemize}
\item $R_\varepsilon$ est positivement transversal \`a $S$, tangent
\`a $T^2 \times \{ a,b\}$;
\item $R_\varepsilon$ vaut $\cos t  {\bdry\over\bdry x} -\sin t
{\bdry\over \bdry y}$
sur $T^2 \times [a+c ,b-c']$ ($\alpha_\varepsilon$ vaut $\cos tdx-\sin tdy$);
\item sur $T^2 \times [a,a+c ]$ et $T^2 \times [b-c' ,b]$,
$R_\varepsilon$ est positivement transversal \`a un feuilletage en anneaux 
$(\F_a \times [a,a+c ])\cup (-\F_b \times [b-c' ,b])$, o\`u $\F_a$
et $\F_b$ sont des
 feuilletages en cercles de $T^2 \times \{ a\}$ et $T^2 \times \{ b\}$
dont les pentes sont $|\varepsilon|$--proches de celles d'une composante
de, respectivement, $\partial S \cap (T^2 \times \{ a\})$ et $\partial 
S\cap (T^2 \times \{ b\})$.
\end{itemize}
\end{theoreme}    

\begin{proof}[D\'emonstration] 
On reprend la d\'emonstration du corollaire~\ref{corollaire : toroidal} 
en construisant sur $V'$ une forme de contact qui trace au bord un
feuilletage caract\'eristique de pente  sup\'erieure 
\`a $\varepsilon$. 
On obtient alors $c$ et $c'$ pour que l'application du lemme~
\ref{lemme : obturation1} donne la transversalit\'e \`a un feuilletage 
en anneaux de pente $|\varepsilon|$--proche de celle de $\partial S$.
\end{proof}

\subsection{Obturations de Dehn sous contr\^ole}

Soit $p/q\in \Q$ et n\'egatif. On r\'ealise une obturation de Dehn de pente
$p/q$ sur une composante $T$ de $\partial V$ pour
obtenir une vari\'et\'e  close $V'$.

\begin{lemme}\label{lemme : obturation2}
Si $|p/q|$ est assez petit, il existe une forme de contact $\beta$
sur $V'$ dont le champ de Reeb est tangent \`a $T$, hypertendu
en restriction \`a $V$ (et transversal \`a une surface
incompressible dans  $V$ qui rencontre $T$) et transversal aux disques
m\'eridiens du tore d'obturation.
\end{lemme}            

\begin{proof}[D\'emonstration]
Il existe une forme de contact hypertendue $\alpha_0$
dont le champ de Reeb $R_0$ est transversal \`a $S$ et tangent \`a $T$.
 Comme $|p/q|$ est assez
petit, on peut choisir un diff\'eomorphisme
de recollement comme suit:

\be              
\item On oriente le cercle m\'eridien $m$ du tore solide
de sorte 
$m$ et $\bdry S$ s'inter\-sectent positivement et transversalement sur
$T$ 
(on rappelle que $T$ est orient\'e comme bord de 
$V$).
\item Le champ de Reeb $R_0$ et $m$ s'intersectent alors n\'egativement
et transversalement sur $T$.
\ee

D'apr\`es le corollaire \ref{corollaire : controle}, il existe
$\varepsilon \in \Q$,  $p/q < \varepsilon <0$, 
pour lequel on peut
trouver une forme de contact $\alpha_\varepsilon$
qui trace sur $T$ un feuilletage de pente $\varepsilon$.
Le champ de Reeb $R_\varepsilon$ de $\alpha_\varepsilon$ 
est de plus colin\'eaire \`a $R_0$. 
Le feuilletage induit sur le tore solide coll\'e
\`a $T$ est un feuilletage en cercles de pente $<0$. 

On choisit par ailleurs un syst\`eme de coordonn\'ees
$D^2 \times \R /\Z =((r,\theta ),z)$ de sorte que $R_\varepsilon$
soit n\'egativement transversal \`a ${\bdry\over \bdry z}$ le long de $T$.

La forme de contact associ\'ee \`a $R_\varepsilon$
vaut $udz+vd\theta$ le long de $T$.

On d\'efinit deux fonctions $f,g \co [0,1]\rightarrow \R^+$
de la m\^eme fa\c con que dans le lemme~\ref{lemme : obturation1}
avec $f(1)=u$, $g(1)=v$, $g(0)=0$, $(f(r),g(r))\neq (0,0)$, $fg'-gf' >0$.

On \'etend alors la forme $\alpha_\varepsilon$ par $f(r)dr+g(r)d\theta$.
Le champ de Reeb associ\'e est transversal aux disques m\'eridiens
du tore solide.
\end{proof}

\section{Cas des fibrations sur le cercle}

Pour rendre cette partie ind\'ependante de
ce qui pr\'ec\`ede, et d\'emontrer
le th\'eor\-\`eme~\ref{theoreme : obturation}, on
n'utilise pas directement la construction du
th\'eor\`eme~\ref{theoreme : existence}, mais une construction
un peu plus simple et explicite, due \`a E~Giroux~\cite{Gi2}.

Soit $S$ une surface compacte, orientable, de bord
non vide, et diff\'erente du disque et de l'anneau.
On note $\phi$ un diff\'eomorphisme de $S$
qui est l'identit\'e pr\`es du bord.
Pour simplifier les notations, on suppose,
comme dans~\cite{Ro1}, que $S$ poss\`ede une seule
composante de bord. La preuve est identique
lorsqu'il y en a plusieurs.

On note $V$ la suspension de $\phi$,
obtenue \`a partir de $P=S\times [0,1]$
en identifiant $(x,1)$ avec $(\phi (x),0)$.
On a une fibration $\pi \co  V\rightarrow S^1$
de fibre $S$.
On note $d\theta$ la forme de longueur sur $S^1$.

Soit \`a pr\'esent $\beta$ une $1$--forme sur
$S$, transversale \`a $\partial S$ et telle
que $d\beta$ est une forme d'aire pour $S$.
Si $\beta_t =(1-t)\beta +t\phi^* \beta$,
$t\in [0,1]$, alors pour $s$ assez petit,
la forme $\alpha=dt +s \beta_t$ est une forme de
contact sur $S\times [0,1]$,
qui induit une forme de contact sur $V$.
Son champ de Reeb $R$ est transversal
aux fibres de $\pi$.

On red\'efinit 
$P=S\times[0,1]$ et $\phi$ afin que la direction verticale
${\bdry\over \bdry s}$, $s\in[0,1]$, soit  tangente
au champ de Reeb $R$, et $V=(S\times[0,1])/
(x,1)\sim(\phi(x),0)$.

\begin{proposition}  \label{roberts}
Il existe une surface branch\'ee $\B$, un voisinage fibr\'e $N(\B )$ 
de $\B$ et $a >0$ tels que:
\be 
\item $N(\B )\subset V$ soit feuillet\'e
verticalement par les orbites de $R$.
\item $N(\B)$ porte, pour tout $\varepsilon
\in [-a,a]$, une lamination $L_\varepsilon$ dont le bord
est feuillet\'e en cercles de pente $\varepsilon$.
\ee
\end{proposition}

\begin{proof}
La preuve est une modification de celle du th\'eor\`eme~4.1 de \cite{Ro1}.   
 Soit $\alpha$ un arc simple, non s\'eparant,
proprement plong\'e dans $S$.
On oriente $\alpha$ et on prend l'orientation induite sur
$\phi(\alpha)$.  Une l\'eg\`ere modification de l'argument
de R.\ Roberts est n\'ec\'essaire car  
$\phi(\alpha)$ et $\alpha$ peuvent ne pas s'intersecter de mani\`ere
{\it efficiente} sur $S$ (i.e.,
le cardinal de $\phi(\alpha)\cap \alpha$ peut ne pas \^etre 
\'egal \`a leur intersection g\'eom\'etrique relative au bord).
Sans perte de g\'en\'eralit\'e, on suppose que
$\phi(\alpha)\pitchfork
 \alpha$.

Il est facile de voir qu'il existe une suite
d'arcs orient\'es, proprement plong\'es et non
s\'eparants
$\alpha_0, \dots,\alpha_m$ sur $S$ avec les propri\'et\'es suivantes:
\be
\item $\alpha_0=\phi(\alpha)$.
\item $\alpha_m$ est isotope \`a $\alpha$ relativement \`a ses extr\'emit\'es.  
\item $\alpha_i$ et $\alpha_{i+1}$ (pour tout $i=0,\dots,m-1$) 
ne s'intersectent pas, ne sont pas parall\`eles et les quatre points
suivants: d\'ebut de $\alpha_{i}$, d\'ebut de
$\alpha_{i+1}$, fin de $\alpha_{i}$, fin de $\alpha_{i+1}$,
apparaissent dans cet ordre  sur $\bdry S$
muni de l'orientation induite.
\ee
Une telle suite est appel\'ee 
``bonne suite positive"
(``{\it positive good sequence}'') dans 
\cite{Ro1}.  Si les quatre points de (3) apparaissent dans
l'ordre inverse sur  
$\bdry S$, la suite est appel\'ee ``bonne suite n\'egative''. 
Il est important de noter qu'on distingue
avec soin un arc et sa classe d'isotopie.

On affirme qu'il existe une extension $\alpha_{m+1},\dots,\alpha_n$
de la bonne suite positive avec $\alpha_n=\alpha$.  
Si on suppose que $\alpha_m$ et $\alpha$ s'intersectent non trivialement
(sauf aux extr\'emit\'es), on utilise un argument ``du plus int\'erieur"
pour trouver des sous-arcs $\beta_m\subset \alpha_m$ et $\beta^-\subset \alpha$ 
qui cobordent un disque $D$ dans $S$ qui n'a pas d'autre intersection
avec $\alpha_m$ et $\alpha_-$ dans son int\'erieur.  Soit $\alpha_{m+4}$ 
l'arc obtenu \`a partir de $\alpha_m$ par isotopie \`a travers $D$ et
soit $\alpha_{m+2}$, l'arc $\alpha_{m+4}$ avec l'orientation oppos\'ee.
Il est ais\'e de trouver $\alpha_{m+1}$ et $\alpha_{m+3}$ (qui
est $\alpha_{m+1}$ avec l'orientation inverse) qui n'intersectent pas 
$\alpha_{m}$, $\alpha_{m+2}$ et $\alpha_{m+4}$. Donc les
paires cons\'ecutives satisfont (3) dans les conditions
d'une bonne suite positive.  L'extension est alors obtenue par induction.

En utilisant une bonne suite positive, on construit \`a pr\'esent $\B$ de
sorte que son voisinage fibr\'e $N(\B)$ soit feuillet\'e par les
orbites du champ de Reeb $R$.  Pour ce faire, on prend $S\times[0,1]$
pour lequel on suppose que $R={\bdry \over \bdry s}$, et
on consid\`ere l'union
$$\B'=\bigcup_{i=0}^n \left(S\times \left\{ {i\over n}\right\}\right)
\cup \bigcup_{i=0}^{n-1} \left(\alpha_i\times\left[{i\over n},
{i+1\over n}\right] \right).$$                 
On munit $S\times \{{i\over n}\}$ de l'orientation de $S$ et
$\alpha_i\times [{i\over n},{i+1\over n}]$ de l'orientation
pour laquelle l'orientation de bord de $\alpha_i\times\{{i\over n}\}$
co\"\i ncide avec l'orientation de $\alpha_i$.
On modifie alors  $\B'$ pr\`es de  $\bdry (\alpha_i\times [{i\over n},{i+1
\over n}])$ pour obtenir une surface branch\'ee
 {\em orient\'ee} $\B$ 
transversale \`a $R$. Pour ce faire, on modifie
$\alpha_i \times [{i\over n},{i+1\over n}]$ par une petite isotopie
pour le rendre positivement transversal \`a $\frac{\partial}{\partial s}$,
puis on d\'eforme l'anneau obtenu au voisinage de son bord
pour introduire deux lignes de branchement avec
$S\times \{ {i\over n}\}$ et $S\times \{ {i+1\over n}\}$.
 On prend alors  pour $N(\B)$ un voisinage
fibr\'e par des arcs de Reeb.
D'apr\`es \cite{Ro1}, il existe une constante positive
$a>0$
telle que, pour tout $\varepsilon\in[0,a)$, le voisinage de 
surface branch\'ee $N(\B)$ porte pleinement une lamination
qui trace un feuilletage en cercles sur $\bdry V$ qui
a pour pente $\varepsilon$.  De fa\c con similaire,
en utilisant une bonne suite n\'egative, il existe
$\B$ et $b<0$ tels que pour tout $\varepsilon\in (b,0]$, $N(\B)$
porte pleinement une lamination dont le bord
est une lamination en cercles de pente $\varepsilon$.
\end{proof}

D'apr\`es la proposition~\ref{roberts},
la lamination $L_\varepsilon$ est transversale \`a $R$.
Le compl\'e\-mentaire de $N(\B )$ dans
$V$ est le produit d'une surface
par une direction verticale tangente \`a $R$.
Toujours d'apr\`es Roberts, la  lamination $L_\varepsilon$
s'\'etend en un feuilletage tendu $\F^V_\varepsilon$, dont la restriction
\`a $N(\B )$ est transversale \`a $R$
et  qui est un feuilletage produit sur $V\setminus N(\B)$,
c'est-\`a-dire qu'il peut \^etre \'egalement pris,
apr\`es isotopie si n\'ec\'essaire, transversal
\` a $R$.
Par construction, $\F^V_\varepsilon \cap \partial V$
est un feuilletage d'holonomie triviale.

Si $\varepsilon$ est assez petit et rationnel,
le feuilletage $\F^V_\varepsilon$ s'\'etend par des disques en un feuilletage
tendu $\F_\varepsilon$ de $V_\varepsilon$, o\`u $V_\varepsilon$
est d\'efinie dans le th\'eor\`eme~\ref{theoreme : obturation}.

Le corollaire~\ref{corollaire : controle} et le
lemme~\ref{lemme : obturation2} permettent de trouver,
lorsque $\varepsilon <0$ (resp. $\varepsilon >0$),
une forme de contact positive (resp. n\'egative) sur $V_\varepsilon$ dont le champ
de Reeb est transversal \`a $\F_\varepsilon$.
Elle est en particulier hypertendue.

Ici, les cons\'equences topologiques dues \`a la pr\'esence
d'une structure hypertendue (ir\-r\'eductibilit\'e, non conjugaison
\`a $S^3$), sont \'egalement impliqu\'ees par la pr\'esence
d'un feuilletage tendu (le rev\^etement universel
de $V_\varepsilon$ est $\R^3$).

La forme de contact est ais\'ee \`a construire. La difficult\'e
est de montrer qu'elle est hypertendue. L'espoir
est d'arriver \`a construire des formes hypertendues
sur une classe de vari\'et\'es plus large, et de fa\c con
plus ais\'ee, que celles qui poss\`edent un feuilletage tendu.

Comme dans \cite{Ro1}, on peut construire
des formes de contact hypertendues, dont le champ de Reeb est transversal
\`a un feuilletage tendu, moyennant des informations
sur les deux premi\`eres surfaces de d\'ecoupages
d'une d\'ecomposition sutur\'ee de $V$.

\begin{theoreme} 
Soit $V$ une vari\'et\'e compacte orientable
bord\'ee par un unique tore. On suppose qu'il existe
un d\'ecoupage de vari\'et\'es sutur\'ees:
$$(V,\partial V) \stackrel {S}\rightsquigarrow (V',\gamma' )
\stackrel {R}\rightsquigarrow
(V'',\gamma''),$$
o\`u $S$ et $R$ sont orientables, $S$ et $R$ ont une seule composante de
bord, $\partial R$ intersecte $S$ en exactement deux arcs non s\'eparants.

Si $(V'',\gamma'')$ est tendue, alors pour $\varepsilon$
assez petit, toute obturation de pente $\varepsilon$ de $V$
porte une structure de contact hypertendue (positive ou n\'egative)
dont le champ de Reeb est transversal \`a un feuilletage
tendu.
\end{theoreme}

La preuve de ce dernier th\'eor\`eme n\'ec\'essite cette fois
d'utiliser la construction du th\'eor\`eme~\ref{theoreme : existence}
dans toute sa g\'en\'eralit\'e.

\section{Annexe}

Dans cette partie, on revient sur l'\'etude du champ
de Reeb construit sur $V_\varepsilon$, dans le cas o\`u
$V$ fibre sur le cercle comme dans la section
pr\'ec\'edente. On reprend int\'egralement les notations
de cette section; en particulier $R_\varepsilon$ est le champ de
Reeb construit dans la section 10.

On peut  montrer directement le th\'eor\`eme suivant, sans recours aux feuilletages.

\begin{theoreme}\label{theoreme : obturation2}
Pour tout $N\in \N$, si $\rho$ est assez petit, pour tout $\varepsilon =(\varepsilon_1 ,...,\varepsilon_n )$
tel que chaque $\varepsilon_i$ soit strictement compris entre $0$ et $\rho$,
aucun rev\^etement de degr\'e inf\'erieur \`a $N$ de $V_\varepsilon$ 
n'est diff\'eomorphe \`a  $S^3$.
De plus, si $V_\varepsilon$ est une sph\`ere d'homologie rationnelle,
elle porte une structure de contact tendue.
\end{theoreme}

On va donner une preuve de ce th\'eor\`eme dans
le cas ou $\partial S$ est connexe, c'est-\`a-dire que
si $\partial S$ est connexe et $\rho$ assez petit,
alors pour tout $\vert \varepsilon \vert <\rho$,
$V_\varepsilon$ n'est pas $S^3$ et porte une structure de contact
tendue dans le cas ou elle est une sph\`ere d'homologie. (La preuve pour les  rev\^etements
de $V_\varepsilon$
n\'ecessite de consid\'erer le cas o\`u $\partial S$
n'est pas connexe.)

En g\'en\'eral, il s'agit de voir que le champ de Reeb $R_\varepsilon$
 ne
poss\`ede pas d'orbite p\'eriodique dans $V_\varepsilon$ qui borde
un disque plong\'e $D$ d'int\'erieur transversal \`a $R_\varepsilon$.
Pour un $N$ fix\'e, cette propri\'et\'e doit aussi
\^etre v\'erifi\'ee, pour $\rho$ assez petit,
sur tout rev\^etement de $V_\varepsilon$ de degr\'e
inf\'erieur \`a $N$. Un tel rev\^etement de $V_\varepsilon$
est \'egalement obtenu par obturation de Dehn
d'une vari\'et\'e qui fibre sur le cercle.
Il suffit donc savoir faire la d\'emonstration dans ce
cas, i.e.\ pour $V_\varepsilon$.

La cl\'e est qu'un tel disque $D$ doit \^etre positivement transversal
\`a l'\^ame des tores d'obturations, qui sont par construction
des orbites de $R_\varepsilon$ (la construction
fournit \'egalement qu'un voisinage de l'\^ame du tore
d'obturation est un produit du disque par des orbites
p\'eriodiques: on peut s'arranger pour que $\partial D\subset V$).

Dans toute la suite, on  suppose, pour simplifier,
que $\partial S$ est connexe. Dans ce cas, on peut supposer
\'egalement que $S$ n'est pas planaire, i.e.\ n'est
pas un disque. La preuve du cas g\'en\'eral suit le m\^eme
sch\'ema, m\^eme si elle comporte  quelques difficult\'es
techniques suppl\'ementaires.

Si on regarde $V\cap D$,
la non-existence de $D$ r\'esulte directement de la proposition suivante
(valable \'egalement lorsque $\partial S$ n'est pas connexe):

\begin{proposition}\label{prop : planaire}
Si $\rho >0$ est assez petit, aucune orbite p\'eriodique $O$ de 
$R$ ne borde de surface planaire $P$ plong\'ee dans $V$, 
$\partial P \setminus O \subset \partial V$, 
dont toutes les composantes de bord orient\'ees $\partial P \setminus O$ 
sont de pente positive et inf\'erieure \`a $\rho$ dans $\partial V$.
\end{proposition}

\begin{proof}[D\'emonstration]
On note $S$ une fibre de $V$, image de $S\times \{ 0\}$ dans la construction.  
Soit $x_0$ un point de $\bdry S$ et $\eta\co [0,1]\rightarrow S$ 
une courbe orient\'ee qui param\'etrise $\bdry S$ et 
satisfait $\eta(0)=\eta(1)=x_0$.

On suppose  que la surface planaire $P$ existe.  
Par une isotopie de $P$, on se ram\`ene au cas 
o\`u toutes les intersections de $\partial P \setminus O$ et de $\partial S$ sont
positives, et o\`u $P$ est transversal \`a $S$.

Comme toutes les intersections de $O$ avec $S$ sont positives,
ce nombre de points d'intersection est aussi l'intersection de $[O ]\in 
H_1 (V;\R )$ et de $[S] \in H_2 (V,\partial V; \R )$. Il est donc \'egalement
le nombre d'intersections de $\partial P \setminus O$
avec $\partial S$.

Ainsi, les composantes \`a bord de $S\cap P$ sont  $p$ arcs
$b_1 ,...,b_p$ allant chacun d'une composante de bord de
$\partial P\setminus O$
\`a $O$.
\`A la source, les pr\'eimages de ces arcs d\'ecoupent la pr\'eimage
de $P$ en un disque dont l'image dans $V$ est not\'ee $D'$.

Toutes les composantes sans bord de $P\cap S$ sont contractibles dans
$S$. En effet, une telle composante $C$
d\'ecoupe dans $P$ une surface planaire
$P_1$ dont le bord ne rencontre pas $O$. L'intersection
homologique de $\partial P_1$ avec $S$ est $0$, l'intersection
de $C$ avec $S$ est $0$, donc $P_1$ n'a pas d'autre composante
de bord que $C$ (toutes les composantes de $\partial P\setminus O$
ont une intersection $+1$ avec $S$).
La courbe $C$ borde un disque dans $P$. Elle est donc contractible
dans $V$  et donc aussi dans $S$ (car $S$ est $\pi_1$--inject\'ee
dans $V$).

On en d\'eduit que la courbe $\partial D'$ est contractible dans
$S\times [0,1]$. Sa projection $\gamma$ dans $S$, par la projection
$\pi$ le long du facteur $[0,1]$, est alors aussi contractible dans $S$.

La courbe $\gamma$ est une concat\'enation d'arcs ferm\'es
\begin{equation}\label{concat}
\gamma =(a_1 b_1 o_1 c_1)(a_2 b_2 o_2 c_2) ... (a_p b_p o_p c_p),
\end{equation}
qui v\'erifient les propri\'et\'es suivantes:

\begin{enumerate}
\item[(i)] Pour tout $1\leq i\leq p$, $a_i \subset \partial S$ est la projection d'une
portion de $\partial P\setminus O$ spiralant dans $\partial S\times [0,1]$ entre $\partial S\times \{0\}$ et $\partial S\times \{1\}$.  Si on ferme les deux extr\'emit\'es de $a_i$ par des arcs $\subset \partial S$ de longueur inf\'erieure (pour une m\'etrique fix\'ee quelconque) \`a la longueur de $\eta$, alors on obtient $A_i$ qui satisfait $[A_i]=[\eta]^{n} \in \pi_1 (S,x_0)$, avec $n \geq 1/\rho -1$.
\item[(ii)] Il existe une permutation $\sigma$ des $p$ premiers entiers $>0$ telle que, 
pour tout $1\leq i\leq p$,  $c_i =\phi (b_{\sigma (i)}^{-1})$.  
Ici $\tau^{-1}$ d\'esigne  l'arc $\tau$ parcouru en sens inverse.
\item[(iii)] Pour tout $1\leq i\leq p$, l'arc $o_i$ est la projection d'une portion de $O$ par $\pi$.
\item[(iv)] Les arcs de la famille $(b_i )_{1\leq i\leq p}$ sont deux \`a deux disjoints, de m\^eme que ceux de la famille $(c_i )_{1\leq i\leq p}$.
\end{enumerate}
La proposition~\ref{prop : homotopie} suivante montre que si $\rho$ est assez petit, une telle courbe $\gamma$ ne peut \^etre homotope \`a z\'ero, ce qui donne une contradiction.
\end{proof}

Dans tout ce qui suit, la surface $S$ est munie d'une m\'etrique hyperbolique, le bord \'etant g\'eod\'esique.  Le diff\'eomorphisme $\phi$ (et pas seulement sa classe d'isotopie) est fix\'e.

\begin{proposition}\label{prop : homotopie}
Soient $\rho, N >0$, et $\gamma$ un lacet de $S$ qui satisfait l'equation \ref{concat} et les conditions {\rm(i)--(iv)} ci-dessus.  Alors, si $\rho$ assez petit, un tel lacet $\gamma$ n'est pas contractible dans $S$.
\end{proposition}

\begin{proof}[D\'emonstration]
L'id\'ee de la preuve est de comptabiliser le nombre de fois que
$\gamma$ ``tourne'' autour de $\partial S$, et de montrer
que, si $\rho$ est assez petit, ce nombre est positif.
On voit dans ce cas que  $\gamma$ ne peut pas \^etre homotope \`a z\'ero.

Pour ce faire, on remarque que chaque arc $a_i$ ``tourne'' au moins $1/\rho -1$ fois autour de $\bdry S$.
Les arcs $o_i$, eux, sont de longueur born\'ee et font donc
autour de $\partial S$ un nombre de ``tours'' born\'e {\it a priori}.
Restent les arcs $b_i$ et les arcs $c_i$ qui sont \'echang\'es
par la monodromie $\phi$: leurs contributions se compensent deux \`a deux
(corollaire~11.5).

La difficult\'e est de pr\'eciser ce que ``tourner'' veut dire.
Pour cela, on ferme d'abord  les arcs $a_i$, $b_i$, $c_i$ et $o_i$
par des arcs de longueur born\'ee {\it a priori}, pour consid\'erer
non plus des arcs mais des lacets bas\'es en $x_0$.  
On note qu'il existe $C>0$ tel que deux points quelconques 
de $S$ peuvent \^etre joints par un segment de longueur inf\'erieure \`a $C$. En ajoutant 
des arcs qui relient des extr\'emit\'es de $a_i$, $b_i$, $c_i$ et $o_i$ \`a $x_0$, 
on obtient des lacets $A_i$, $B_i$, $O_i$ et $C_i$.  
Ici, $A_i$ est d\'efini comme dans la condition (i).  
Alors $\gamma$ est homotope \`a un lacet:
\begin{equation}\label{gammanul}
\gamma_0=(A_1B_1 O_1C_1)(A_2B_2O_2C_2)\dots (A_p B_p O_p C_p),
\end{equation}
o\`u $[A_i]=[\eta]^n$.

Le nombre de ``tours'' fait par un lacet $L$ autour $\partial S$ 
se traduit alors par l'apparition de puissances de la 
classe $[\eta] \in  \pi_1 (S,x_0)$ dans l'\'ecriture 
r\'eduite de $[L]\in \pi_1 (S,x_0)$ pour un certain syst\`eme 
de g\'en\'erateurs de $\pi_1 (S,x_0)$.  On remarque ici que $\pi_1(S,x_0)$ est un groupe libre.
On choisit un syst\`eme $\langle \alpha_1,\dots,\alpha_{2g(S)} \rangle$ 
de g\'en\'erateurs de $\pi_1(S,x_0)$ comme suit.
On se donne un arc {\em orient\'e} g\'eod\'esique $g$, d'extr\'emit\'es incluses dans $\bdry S$ et
non homotope \`a un arc de $\partial S$.  
On suppose sans perte de g\'en\'eralit\'e
(quitte \`a modifier $\phi$) que $\phi (g)\pitchfork g$.  
On prend un syst\`eme tel qu'un seul g\'en\'erateur $\alpha_1$ intersecte $g$, et 
que cette intersection soit  r\'eduite \`a un point.  
On souligne le statut particulier de $[\eta]$, qui se 
d\'ecompose (sans perte de g\'en\'eralit\'e) de la fa\c con suivante:
\begin{equation}\label{commutateur}
\eta\sim [\alpha_1,\alpha_2]\cdots [\alpha_{2g(S)-1},\alpha_{2g(S)}].
\end{equation}
 Ici, $\sim$ est la relation d'homotopie.

Il faut alors montrer que, pour $\rho$ assez petit, l'\'ecriture
r\'eduite de $[\gamma_0]\in \pi_1 (S,x_0)$  dans ce syst\`eme
de g\'en\'erateurs comprend un nombre strictement positif de facteurs
$[\eta]$, ce qui montre que $[\gamma_0]\neq 1$.

Soit $[B] \in \pi_1 (S,x_0 )$. On s'int\'eresse \`a l'occurence des facteurs $[\eta]$ dans l'\'ecriture r\'eduite de $[B]$ dans ce syst\`eme de g\'en\'erateurs.  Si l'\'ecriture de $[B]$ n'est pas r\'eduite, on compte les facteurs de la m\^eme mani\`ere, avec la convention qu'il peut y avoir des facteurs d'exposant $0$. Par exemple, $[\eta]^2 [\eta]^{-1} [\eta ]^{-1}$ sera d\'ecrit comme $[\eta ]^0$ et comptabilis\'e comme un groupement.

On note $\tilde{S}$ le rev\^etement universel de $S$, et $\tilde{\eta}$ (resp. $\tilde{g}$)
le rappel de $\eta$ (resp. $g$) dans $\tilde{S}$. Pour chaque composante $\tilde{d}$ de $ {\tilde{\eta}}$, on consid\`ere l'ensemble $\tilde{g}_{\tilde{d}}$ des composantes de
$\tilde{g}$ qui rencontrent $\tilde{d}$ 
et dont l'orientation au point d'intersection pointe dans $S$.

\begin{lemme}\label{lemme : puissance}
Soit $[B] \in \pi_1 (S,x_0 )$. Un mot $[\eta ]^m$, $m\in \Z$ et $|m|$ maximal,
appara\^\i t dans l'\'ecriture r\'eduite de $[B]$ si et seulement
si le rappel $\tilde{B}$ de $B$ dans $\tilde{S}$ a une intersection
homologique de $m-1$, $m$ ou $m+1$ avec l'un des $\tilde{g}_{\tilde{d}}$.
\end{lemme}

\begin{proof}[D\'emonstration]
 D'abord, on homotope $B$ \`a son \'ecriture sous forme r\'eduite.  
L'intersection homologique avec $\tilde{g}_{\tilde{d}}$ reste 
invariante pendant l'homotopie. On garde la notation $B$ pour ce nouveau lacet.  
Sans perte de g\'en\'eralit\'e, on suppose que $m>1$.

Ensuite, on observe que les seuls rappels des repr\'esentants 
des g\'en\'erateurs qui rencontrent $\tilde{g}_{\tilde{d}}$ 
sont des relev\'es de $\alpha_1$.  On suppose qu'il existe 
un relev\'e $\tilde \alpha_1$ de $\alpha_1$ dans le rappel 
$\tilde{B}$ de $B$ qui intersecte $\tilde{g}_{\tilde{d}}$.  
Alors $\tilde \alpha_1^{-1}$ n'appara\^\i t pas dans $\tilde{B}$, 
car $B$ est sous forme r\'eduite et $\tilde S$ est contractible.  
L'ensemble $\tilde{g}_{\tilde{d}}$ d\'ecoupe $\tilde S$ en 
composantes adjacentes $\{S_i\}_{i\in\Z}$, o\`u $S_i$ et $S_{i+1}$ 
bordent la composante $\tilde{g}_i$ de $\tilde{g}_{\tilde d}$. 
(Ici, l'indice $i$ cro\^\i t dans la direction de l'orientation de $\bdry \tilde S$.) 
Alors, si $\tilde{B}$ a une intersection homologique 
$m$ avec $\tilde{g}_{\tilde{d}}$, 
le point initial de $\tilde B$ est (sans perte de g\'en\'eralit\'e) 
dans $S_0$ et le point terminal de $\tilde B$ est dans $S_m$. 
La premi\`ere et la derni\`ere intersection de $\tilde{g}_{\tilde{d}}$ 
avec $\tilde{B}$ sont $\tilde{g}_0\cap \tilde\alpha_1^1$ et $\tilde{g}_{m-1}\cap \tilde\alpha_1^2$, 
o\`u $\tilde\alpha_1^1$ et $\tilde\alpha_1^2$ sont les  
relev\'es de $\alpha_1$ qui partent de $\tilde{d}$.  
Les arcs $\tilde\alpha_1^1$ et $\tilde\alpha_1^2$ d\'elimitent $\tilde{B_0}$, qui est le 
relev\'e d'un sous-mot $B_0$ dans l'\'ecriture r\'eduite de $B$. 
(Ici, les arcs $\tilde\alpha_1^1$ et $\tilde\alpha_1^2$ sont non inclus dans $\tilde{B_0}$.)  
Comme $\tilde{S}$ est contractible, $\alpha_1 B_0=\eta^{m-1}$.
 Le mot $B=\dots (\alpha_1 B_0) \alpha_1\dots= \dots \eta^{m-1} \alpha_1\dots$ 
appara\^\i t dans l'\'ecriture r\'eduite, gr\^ace \`a l'equation~\ref{commutateur}.

Si l'intersection homologique vaut $m$, alors $\eta^{m-1}$ 
appara\^\i t dans l'\'ecriture\break r\'eduite de $B$. 
R\'ecipro\-que\-ment, si un sous-mot $\eta ^m$ appara\^\i t 
dans l'\'ecriture r\'eduite de $B$, alors $\eta^m$ contribue pour $m$ 
\`a l'intersection homologique.  
Gr\^ace au paragraphe pr\'ec\'edent, l'intersection 
homologique entre $\tilde B$ et $\tilde{g}_{\tilde d}$ est  $\geq m$.
\end{proof}

La caract\'erisation pr\'ec\'edente permet de montrer 
que l'apparition d'un facteur $[\eta]^m$ dans 
l'\'ecriture r\'eduite d'un \'el\'ement $[B] \in \pi_1 (S,x_0 )$
est une propri\'et\'e invariante, \`a une puissance 
inf\'erieure \`a un certain $r_\phi$ pr\`es, par le diff\'eomorphisme $\phi$.

\begin{corollaire}\label{corollaire : isometrie}
Il existe $r_\phi >0$ tel qu'une puissance $[\eta ]^m$ appara\^\i t 
dans l'\'ecriture r\'eduite d'un mot $[B] \in \pi_1 (S,x_0 )$ 
si et seulement si une puissance $[\eta ]^{m+k}$, $-r_\phi \leq k\leq r_\phi$, appara\^\i t dans l'\'ecriture r\'eduite de $[\phi (B)]$.
\end{corollaire}

\begin{proof}[D\'emonstration]
Les rappels de $g$ et de $\phi (g )$ dans $\tilde{S}$
se rencontrent en un nombre fini de points, born\'e par
$M$ ind\'ependemment de $\eta$ et des rappels consid\'er\'es.
On note $\tilde{\phi}$ le relev\'e de $\phi$ \`a $\tilde{S}$,
qui est l'identit\'e sur $\tilde{d}$.

L'intersection de $\tilde{\phi} (\tilde{B} )$ avec
$\tilde{\phi} (\tilde{g}_{\tilde{d}} )$ en homologie, qui est la
m\^eme que celle de $\tilde{B}$ avec
$\tilde{g}_{\tilde{d}}$ est donc \'egalement
la m\^eme que celle de $\tilde{\phi} (\tilde{B} )$
avec $\tilde{g}_{\tilde{d}}$, corrig\'ee d'un facteur
$k\in [-M-2,M+2]$.

Par application du lemme~\ref{lemme : puissance},
les sous-mots $[\eta]^m$ qui apparaissent dans\break
l'\'ecriture r\'eduite de $[\phi (B)]$ sont donc les m\^emes,
\'eventuellement \`a $M+4$ unit\'es de puissance pr\`es,
que ceux qui apparaissent dans celle de $[B]$.
On peut prendre $r_\phi =M+4$.
\end{proof}

\begin{remarque}
{\rm En particulier, de nombreux sous-mots $[\eta ]^k$, $-r_\phi \leq k\leq r_\phi$ 
peuvent \^etre pr\'esent dans une \'ecriture et pas dans l'autre.}
\end{remarque}

Le lemme qui suit va permettre d'exploiter le fait que les arcs $b_i$ et $c_i$, $1\leq i\leq p$,
sont plong\'es, car $P$ l'est.

\begin{lemme}\label{lemme : plonge}
Soit $\delta$ un arc orient\'e plong\'e dans $S$ et $\tilde\delta$ un relev\'e dans $\tilde{S}$. Il existe au plus deux composantes de bord $\tilde d_1$ et $\tilde d_2$ de $\partial \tilde{S}$ telles que la somme des intersections alg\'ebriques entre $\tilde\delta$ et $\tilde{g}_{\tilde d_i}$ soit, en valeur absolue, sup\'erieure ou \'egale \`a $5$.
\end{lemme}

\begin{proof}[D\'emonstration]
On suppose qu'il existe deux telles composantes $\tilde d_1$ et $\tilde d_2$. 
On montre qu'il n'y en a pas d'autre. Ici, on utilise 
les deux propri\'et\'es suivantes: (i) $\delta$ est {\em plong\'e} 
dans $S$ et (ii) si $\tilde\delta$ intersecte $\tilde{g}_{\tilde d}$ 
plus d'une fois, o\`u $\tilde d$ est un relev\'e de $\bdry S$, 
alors le sous-arc correspondant de $\delta$ spirale parall\`element \`a $\bdry S$.

Les nombres d'intersection consid\'er\'es sont invariants par isotopie de $\delta$ relative
\`a son bord. On a le fait suivant (la d\'emonstration suit un argument 
standard o\`u on cherche le sous-arc le plus int\'erieur):

\s\n
{\bf Fait 1}\qua Il existe une isotopie de $\delta$ relative 
\`a son bord en un arc plong\'e $\delta_0$ telle que 
toutes les intersections de $\tilde\delta_0$, obtenu 
par rel\`evement de l'isotopie \`a $\tilde\delta$, 
avec $\tilde d_1$ et $\tilde d_2$ soient de m\^eme signes.

\s
Soient, pour $i=1,2$, $\tilde p_i$ (resp.\ $\tilde q_i$) 
la premi\`ere (resp.\ la derni\`ere) intersection de $\tilde{g}_{\tilde d_i}$ 
avec $\tilde\delta_0$. On note $p_i$ (resp.\ $q_i$) 
la projection de $\tilde p_i$ (resp.\ $\tilde q_i$) 
dans $S$. Comme $\delta_0$ est plong\'e, les rappels 
de $p_i$ et $q_i$ sur une m\^eme composante de $\tilde{g}_{\tilde d_i}$ 
sont distincts. On suppose, quitte \`a \'echanger les 
points $p_i$ et $q_i$, que celui de $q_i$ est situ\'e entre celui de $p_i$ et $\tilde d_i$.

Pour chaque $i$ on d\'efinit $\delta_1$ comme 
l'adh\'erence de la composante de $\delta_0 \setminus p_i$ qui 
contient $q_i$ et $\tilde\delta_1 \subset \tilde\delta$ son rappel. On note de plus $\tilde\delta_2 \subset \tilde\delta_1$ l'arc bord\'e par $\tilde p_i$ et $\tilde q_i$.

\s\n
{\bf Fait 2}\qua Si $\tilde{d}$ est un relev\'e de $\bdry S$ 
dans $\tilde{S}$, distinct de $\tilde d_1$ et $\tilde d_2$, 
alors l'intersection de $\tilde \delta_1$ dans $\tilde{S}$ 
avec $\tilde{g}_{\tilde d}$ est en valeur absolue inf\'erieure \`a $1$.

\begin{proof}[Preuve du fait 2]
L'arc $\tilde\delta_2$ est isotope, relativement \`a son bord, 
\`a un arc contenu dans $\tilde{g}_{\tilde d_i}\cup \tilde{d_i}$. 
Son intersection avec $\tilde{g}_{\tilde d}$ est $0$ ou $\pm 1$, 
car c'est le cas pour $\tilde d_i$ ($S$ est hyperbolique). 
Ce nombre d'intersection est de plus invariant par isotopie relative au bord.

Par d\'efinition, $\tilde \delta_1 \setminus \tilde\delta_2$ 
ne rencontre pas $\tilde{g}_{\tilde d_i}$. Il est situ\'e 
entre deux composantes successives $\tilde g_1$ et 
$\tilde g_2$ de $\tilde{g}_{\tilde d_i}$. Les points 
$\tilde g_j \cap \tilde d_i$ d\'elimitent un sous-arc 
$\tilde\beta$ dans $\tilde d_i$. Un translat\'e $t(\tilde \delta_1)$ 
de $\tilde\delta_1$ parall\`ele \`a $\tilde{d_i}$ rencontre les deux composantes $\tilde g_1$
et $\tilde g_2$ en deux points $r_1$ et $r_2$.  La courbe lisse 
par morceaux compos\'ee de l'union de $\tilde\beta$, des 
sous-arcs de $\tilde g_1$ et $\tilde g_2$ situ\'es entre $r_1$ et $r_2$ 
et l'intersection avec $\tilde{d_i}$, et le sous arc de $t(\tilde\delta_1)$ 
d\'elimit\'e par $r_1$ et $r_2$ est sans point double. Elle borde un 
disque $D$ dans $\tilde{S}$. Comme $\delta_0$ est plong\'e, 
$\delta_1 \setminus \delta_2 \subset D$. L'intersection de $\beta$ 
avec $\tilde{g}_{\tilde d}$ vaut au plus $1$ en valeur absolue. 
C'est donc \'egalement le cas de l'intersection de 
$\tilde\delta_1 \setminus \tilde\delta_2$ avec $\tilde{g}_{\tilde d}$.

Finalement, l'intersection de $\tilde{g}_{\tilde{d}}$ avec 
$\tilde\delta_1$ est en valeur absolue major\'ee par $2$ 
(en fait m\^eme, par l'intersection de $\tilde{g}_{\tilde d}$ 
avec $\tilde d_i$, c'est-\`a-dire $1$). \end{proof}

Cette d\'emonstration donne \'egalement que 
les intersections de $\tilde{\delta}$ avec 
$\tilde{g}_{\tilde d}\not= \tilde{g}_{\tilde d_i}$, $i=1,2$, sont toutes (\`a une unit\'e pr\`es) contenues dans $\tilde\delta_0 \setminus \tilde\delta_1$ 
et que l'intersection de $\tilde\delta_0 \setminus \tilde\delta_1$ 
avec $\tilde{g}_{\tilde{d}}$ est en valeur absolue major\'ee par $2$.
\end{proof}

La longueur des arcs $O_i$ est inf\'erieure \`a une 
certaine constante (d\'ependant de $\phi$). L'\'ecriture 
r\'eduite de $O_i$ dans le syst\`eme de g\'en\'erateurs 
choisi est donc de longueur born\'ee ind\'ependemment de $\rho$ et $p$.  
De la m\^eme mani\`ere, on obtient l'existence de $\mathfrak{d}\in \N$, 
ind\'ependant de $\rho$ et $p$, tel que pour tout sous-mot $\eta^k$ 
apparaissant dans $A_i=\eta^n$, $B_i$, $C_i$ ou $O_i$, alors $-\mathfrak{d}< k< \mathfrak{d}$, 
sauf \'eventuellement pour au plus $2$ sous-mots dans chaque $B_i$ et 
$C_i$ (lemme~\ref{lemme : plonge}) et exactement $1$ dans chaque $A_i$.

On applique \`a pr\'esent le corollaire~\ref{corollaire : isometrie} 
\`a $B_i$ et $\phi(B_i^{-1})$. L'\'ecriture r\'eduite de $\phi (B_i^{-1})$ fait intervenir 
$\eta$ avec des exposants situ\'es \`a distance born\'ee $r_\phi$ 
(ind\'ependemment de $\rho$ et $p$) 
des oppos\'es de ceux de $B_i$. Quitte \`a changer la valeur 
de $\mathfrak{d}$ de sorte  que $\mathfrak{d}\geq r_\phi$, on suppose, 
sans perte de g\'en\'eralit\'e, qu'au plus $2$ sous-mots dans 
l'\'ecriture de $\phi (B_i^{-1})$ (et de $C_i$) ont un exposant 
en valeur absolu $\geq \mathfrak{d}$.

On calcule la somme totale $C$ des exposants des sous-mots (``maximaux'') 
de type $\eta^k$, o\`u $k\geq N\mathfrak{d}$, dans l'\'ecriture r\'eduite de $\gamma_0$.  
(La constante $N$ reste \`a d\'eterminer.) Un sous-mot maximal 
de type $\eta^k$, $k\geq N\mathfrak{d}$ est dit {\em long}. Soit $c$ 
la somme totale des exposants des sous-mots longs pr\'esents dans les \'ecritures r\'eduites 
de $A_i=\eta^n$, $B_i$, $C_i$ et $O_i$, c'est-\`a-dire 
la somme totale du nombre de sous-mots apparaissant dans l'\'ecriture de $\gamma_0$ avant simplification.  Comme chaque $A_i$ contribue pour  $n$, chaque paire $B_i$ et $C_i$ contribue 
pour au moins $-2\mathfrak{d}$ et chaque $O_i$ contribue pour  $0$, on a $c\geq p(n-2\mathfrak{d})$.

Maintenant, on simplifie l'\'ecriture de $\gamma_0$ 
donn\'e par l'\'equation~\ref{gammanul}. Les 
simplifications s'effectuent n\'ec\'e\-ssairement aux 
jonctions des mots r\'eduits. On commence par simplifier \`a la premi\`ere jonction, puis ensuite \`a la premi\`ere jonction restante....
\`A une jonction:
\be
\item[(a)] soit on fait appara\^\i tre, ou dispara\^\i tre $[\eta]^{\pm 1}$.
\item[(b)] soit un sous-mot $[\eta]^r$ d'un c\^ot\'e s'adjoint avec un sous-mot $[\eta]^s$
de l'autre, et soit ($\mbox{b}_1$) $s\neq -r$ et le processus,
pour cette jonction, s'arr\^ete et donne naissance \`a un sous-mot $[\eta ]^{r+s}$,
soit ($\mbox{b}_2$) $r=-s$, et on \'elimine un sous-mot, la somme totale des exposants restant fixe.
\ee
Les cas (a) et ($\mbox{b}_1$) arrivent au plus une fois chacun \`a chaque jonction.

Si un sous-mot long ne participe pas \`a une jonction, alors $c$ 
ne change pas. Si on prend $N\geq 5$, on ne peut pas cr\'eer 
un sous-mot long par une concat\'enation successive de $\eta^k$ 
avec $k<\mathfrak{d}$. Supposons qu'un ou deux sous-mots longs 
participent \`a une jonction.  Alors, l'erreur 
est major\'ee par $\mathfrak{d}$, et on obtient:
\begin{equation}
C\geq p ((n-2\mathfrak{d})-4\mathfrak{d})=p(n-6\mathfrak{d}).
\end{equation}
Si $n>>6\mathfrak{d}$ o\`u ${1\over \rho}>>6\mathfrak{d}$, on a $C>>0$.  
Donc le lacet $\gamma_0$ n'est pas contractible.
\end{proof}

\end{document}

%% file: gtoutput.tex
\def\ifplaintex{\expandafter\ifx\csname documentclass\endcsname\relax}

\ifplaintex 
\else
\fi

\expandafter\ifx\csname beginpicture\endcsname\relax
\expandafter\ifx\csname documentclass\endcsname\relax
\input pictex \else
\input prepictex \input pictex \input postpictex \fi\fi

\def\gt{{\mathsurround=0pt\it $\cal G\mskip-2mu$eometry \&\ 
$\cal T\!\!$opology}}        

\def\gtp{{\mathsurround=0pt\it $\cal G\mskip-2mu$eometry \&\ 
$\cal T\!\!$opology $\cal P\!$ublications}}  

\def\lognumber#1{\def\thelognumber{#1}}
\def\volumenumber#1{\def\thevolumenumber{#1}}
\def\papernumber#1{\def\thepapernumber{#1}}
\def\volumeyear#1{\def\thevolumeyear{#1}}

\def\pagenumbers#1#2{\def\startpage{#1}\def\finishpage{#2}}
\def\published#1{\def\publishdate{#1}}
\def\proposed#1{\def\theproposer{#1}}
\def\seconded#1{\def\theseconders{#1}}
\def\received#1{\def\receiveddate{#1}}
\def\revised#1{\def\reviseddate{#1}}
\def\accepted#1{\def\accepteddate{#1}}
\def\asciititle#1{\def\theasciititle{#1}}
\def\covertitle#1{\def\thecovertitle{#1}}

\def\asciiaddress#1{\def\theasciiaddress{#1}}
\def\asciiemail#1{\def\theasciiemail{#1}}

\long\def\asciiabstract#1{\long\def\theasciiabstract{#1}}

\let\\\par\let\thelognumber\relax
\let\thevolumenumber\relax\let\thepapernumber\relax
\let\thevolumeyear\relax\let\thesamplenumber\relax\let\startpage\relax
\let\finishpage\relax\let\publishdate\relax\let\receiveddate\relax
\let\reviseddate\relax\let\accepteddate\relax\let\theasciititle\relax
\let\thecovertitle\relax\let\theasciiauthors\relax\let\theasciiaddress\relax
\let\theasciiabstract\relax
\let\theasciiemail\relax\let\theshortauthors\relax\let\theshorttitle\relax

\long\def\maketitlep{   

\count0=\startpage

\gt\hfill      
\beginpicture
\setcoordinatesystem units <0.33truein, 0.33truein> point at 2.2 0.9
\setplotsymbol ({$\cal G$})
\plotsymbolspacing=9truept
\circulararc 315 degrees from 0 1 center at 0 0
\setplotsymbol ({$\cal T$})
\circulararc 315 degrees from 1 -1 center at 1 0
\endpicture
\break
{\small\ifx\thesamplenumber\relax 
Volume \else Sample
\fi\thevolumenumber\ (\thevolumeyear)
\startpage--\finishpage\nl
Published: \publishdate}
\vglue 0.5truein plus 0.4fil minus 0.1truein

{\parskip=0pt\leftskip 0pt plus 1fil\def\\{\par\smallskip}{\ifplaintex\large
\else\Large\fi\bf\thetitle}\par\medskip}   

\vglue 0pt plus 0.1fil 

{\parskip=0pt\leftskip 0pt plus 1fil\def\\{\par}{\sc\theauthors}
\par\medskip}

\vglue 0pt plus 0.1fil 

{\small\parskip=0pt\let\newline\\
{\leftskip 0pt plus 1fil\def\\{\par}{\sl\theaddress}\par}
\expandafter\ifx\theemail\relax    
\relax\else\vglue 5pt plus 0.02fil minus 2pt\def\\{\stdspace{\rm 
and}\stdspace} 
\cl{Email:\stdspace\tt\theemail}\fi
\ifx\theurl\relax                  
\relax\else\vglue 5pt plus 0.02fil minus 2pt\def\\{\stdspace{\rm 
and}\stdspace}
\cl{URL:\stdspace\tt\theurl}\fi\par}

\vglue 7pt plus 0.3fil minus 3pt

{\bf Abstract}
\vglue 5pt plus 0.1fil minus 2pt

\theabstract

\vglue 7pt plus 0.3fil minus 3pt

{\bf AMS Classification numbers}\quad Primary:\quad \theprimaryclass

Secondary:\quad \thesecondaryclass

\vglue 5pt plus 0.3fil minus 2pt

{\bf Keywords:}\quad \thekeywords

\vglue 10pt plus 0.5fil minus 5pt

{\small  Proposed: \theproposer\hfill Received: \receiveddate\nl
Seconded: \theseconders\hfill 
\ifx\reviseddate\relax                         
Accepted: \accepteddate                        
\else
Revised: \reviseddate                          
\fi}
\eject
}       

\let\maketitlepage\maketitlep
\let\maketitle\maketitlepage

\font\phead=cmsl9 scaled 950
\font=cmsl9 scaled 1050
\font\pnum=cmbx10 scaled 913
\font\lnum=cmbx10 
\font\pfoot=cmsl9 scaled 950
\font=cmsl9 scaled 1050
\ifplaintex
\headline{\vbox to 0pt{\vskip -4.5mm\line{\small\phead\ifnum
\count0=\startpage ISSN 1364-0380 (on line)
1465-3060 (printed) \hfill {\pnum\folio}\else\ifodd\count0\def\\{ }%
\ifx\theshorttitle\relax\thetitle\else\theshorttitle\fi\hfill{\pnum\folio}
\else\def\\{ and }{\pnum\folio}\hfill\ifx\theshortauthors\relax\theauthors
\else\theshortauthors\fi\fi\fi}\vss}}
\footline{\vbox to 0pt{\vglue 0mm\line{\small\pfoot\ifnum\count0=\startpage
\copyright\ \gtp\hfill\else
\gt, Volume \thevolumenumber\ (\thevolumeyear)\hfill\fi}\vss
}}
\else
\makeatletter
\def\@oddhead{{\small\lhead\ifnum\count0=\startpage ISSN 1364-0380 (on line)
1465-3060 (printed) \hfill {\lnum\number\count0}\else\ifodd\count0
\def\\{ }\ifx\theshorttitle\relax \thetitle \else\theshorttitle\fi\hfill
{\lnum\number\count0}\else\def\\{ and }{\lnum\number\count0}
\hfill\ifx\theshortauthors\relax 
\theauthors\else\theshortauthors\fi\fi\fi}}\def\@evenhead{\@oddhead}
\def\@oddfoot{\small\lfoot\ifnum\count0=\startpage\copyright\ \gtp\hfill\else
\gt, Volume \thevolumenumber\ (\thevolumeyear)\hfill\fi}
\def\@evenfoot{\@oddfoot}
\makeatother
\fi

\newwrite\gtoutfile
\long\gdef\makeheadfile{  
{\def\\{, }\def\s{ }
\immediate\openout\gtoutfile head.xxx
\immediate\write\gtoutfile{Proxy-for: \ifx\theasciiauthors\relax
\theauthors\else\theasciiauthors\fi\s<\ifx\theasciiemail\relax\theemail\else\theasciiemail\fi>}
\immediate\write\gtoutfile{\noexpand\\}
\immediate\write\gtoutfile{Authors: \ifx\theasciiauthors\relax
\theauthors\else\theasciiauthors\fi}
{\def\\{ }\immediate\write\gtoutfile{Title: \ifx\theasciititle\relax
\thetitle\else\theasciititle\fi}}
\immediate\write\gtoutfile{Subj-class: GT or SG or MG etc}
\immediate\write\gtoutfile{MSC-class: \theprimaryclass\ifx\thesecondaryclass\relax\else, \thesecondaryclass\fi}
\immediate\write\gtoutfile{Journal-ref: Geom. Topol. \thevolumenumber
(\thevolumeyear) \startpage-\finishpage}
\immediate\write\gtoutfile{Comments: Published by Geometry and Topology at}
\immediate\write\gtoutfile{\s\s http://www.maths.warwick.ac.uk/gt/GTVol\thevolumenumber/paper\thepapernumber.abs.html}
\immediate\write\gtoutfile{\noexpand\\}
\immediate\write\gtoutfile{}
\ifx\theasciiabstract\relax
\immediate\write\gtoutfile{\theabstract}\else
\immediate\write\gtoutfile{\theasciiabstract}\fi
\immediate\write\gtoutfile{}
\immediate\write\gtoutfile{\noexpand\\}
\immediate\write\gtoutfile{}
\immediate\closeout\gtoutfile}}  

\def\maketitlepage{\maketitlep\makeheadfile}
\let\maketitle\maketitlepage